\documentclass{article}

\usepackage[english]{babel}

\usepackage[letterpaper,top=2cm,bottom=2cm,left=3cm,right=3cm,marginparwidth=1.75cm]{geometry}

\usepackage{amsmath, amssymb}
\usepackage{amsthm}  
\usepackage{graphicx}
\usepackage{booktabs}
\usepackage{subfigure}
\usepackage{color}
\usepackage[colorlinks=true, allcolors=blue]{hyperref}
\usepackage{hyperref}
\usepackage{cleveref}
\usepackage{bm}

\newtheorem{lemma}{Lemma} 
\newtheorem{definition}{Definition} 
\newtheorem{theorem}{Theorem} 
\newtheorem{remark}{Remark} 

\definecolor{zhu}{rgb}{0.000,0.000,1.000}
\definecolor{yin}{rgb}{1.000,0.000,0.000}

\title{An adaptive finite element method for the interior problem of Helmholtz equations with generalized transmission boundary conditions}
\author{Zhe Zhu\thanks{School of Mathematics and Computational Science, Xiangtan University, Xiangtan, Hunan, 411105, China; School of Engineering and Management, Pingxiang University, Pingxiang, Jiangxi, 337000, China(\texttt{zhuzheedu@163.com})},\;
Lei Lin\thanks{State Key Laboratory of Mathematical Sciences and Institute of Computational Mathematics and Scientific/Engineering Computing, Academy of Mathematics and Systems Science, Chinese Academy of Sciences, Beijing 100190, China (\texttt{linlei@amss.ac.cn})},\;
Wei Yang\thanks{School of Mathematics and Computational Science, Xiangtan University, Xiangtan, Hunan, 411105, China           (\texttt{yangwei@xtu.edu.cn})},\;
Tao Yin\thanks{State Key Laboratory of Mathematical Sciences and Institute of Computational Mathematics and Scientific/Engineering Computing, Academy of Mathematics and Systems Science, Chinese Academy of Sciences, Beijing 100190, China (\texttt{yintao@lsec.cc.ac.cn})}}
\date{\vspace{-5ex}}

\begin{document}
\maketitle
\begin{abstract}
This paper proposes a new adaptive finite element method for solving the interior transmission problem with generalized transmission boundary conditions (GTBCs) involving a surface Laplace-Beltrami operator. A $\mathbb{T}$-coercivity technique relying on a lifting operator is established to show the well-posedness of the variational problem in appropriate nonstandard Sobolev spaces. Then both the a priori error estimates and a posteriori error estimates including the upper bound and local lower bound of the finite element discretization are constructed. Numerical examples are presented to illustrate the accuracy and convergence of the adaptive finite element method. Extension and high efficiency of the proposed method for the computation of the interior transmission eigenvalues for the problem with GTBCs are also discussed numerically.
\end{abstract}
\textbf{Keywords:}
transmission problem; Laplace-Beltrami operator; $\mathbb{T}$-coercivity; adaptive finite element method; error estimates.

\section{Introduction}

The classical transmission problems \cite{CK19} modeled by the wave equations and continuity transmission boundary conditions of the Cauchy data can be used to characterize the wave scattering by penetrable obstacles in many practical applications. These problems, together with the related interior transmission problems (ITPs) arising from the study of reconstruction of obstacles or medium parameters \cite{CCH16,CGH10}, have been extensively studied in open literature for the problems in acoustics \cite{CGH10}, electromagnetics \cite{CN21} and elastodynamics \cite{WBY24}. However, for the problems with thin coatings \cite{BC16} or containing complex surface effects \cite{DHJ12}, the numerical evaluation of the classical transmission problems requires dealing with multiscale geometries, which suffers from expensive computing costs. Hence, some equivalent models with generalized transmission boundary conditions (GTBCs) involving high-order derivatives have been introduced to approximate the local effects and avoid the discretization of thin layers. The same idea has also been developed for the wave scattering by strongly absorbing structures or impenetrable obstacles with thin coatings and then problems with generalized impedance boundary conditions (GIBCs) result \cite{AHS11,AS14,HJN05,HJN08,IYLL21,KLL16}. This paper is interested in the corresponding ITP coming from the study of inverse transmission problems with GTBCs \cite{GHK26} which involves a Laplace-Beltrami operator on the boundary and is devoted to proposing an efficient adaptive finite element method (FEM) with rigorous numerical analysis.

It should be pointed out that the GTBCs considered in this work involve a second-order Laplace-Beltrami differential operator, which gives rise to several challenges in both well-posedness analysis and error estimates. For the well-posedness analysis, on one hand, the standard $H^1$ space used in classical ITPs \cite{BBCH11,WC13} is insufficient to define an appropriate variational formulation for the generalized ITP and then a new nonstandard Sobolev space which can characterize both the bulk- and surface-regularities is necessary. On the other hand, the variational formulation of ITP does not provide compact perturbation of a coercive term to allow a direct application of the Fredholm alternative. Motivated by the study of the generalized Robin boundary value problems~\cite{KCDQ15}, we introduce a nonstandard Sobolev-space setting that incorporates the additional surface regularity induced by the GTBCs and then establish the well-posedness result by means of a suitable $\mathbb T$-coercivity argument. The $\mathbb T$-coercivity technique originally proposed in \cite{BBCZ10} gives an innovative way to reformulate the variational problem using a so-called $\mathbb T$-operator and its effectiveness in the study of ITPs has been shown in~\cite{BBCH11}. To treat the ITP with GTBCs, a lifting operator relating to a well-posed variational problem can be introduced~\cite{CC13} to define the $\mathbb T$-operator. One interesting thing is that owing to the additional surface regularity, the argument in well-posedness analysis, which can also be extended to the anisotropic problems, is valid under
less restrictive assumptions on the parameters compared with the results for the classical ITPs~\cite{NN17,WC13}. 

Turning to the finite element error analysis, the adaptive approximation for problems with non-smooth boundaries, which may significantly reduce the solutions' regularities, involves not only the standard bulk equation residuals but also additional surface residuals with GTBCs and tangential flux contributions induced by the Laplace-Beltrami operator. As far as we know, only a few results on the a posteriori error estimates for the problems with boundary conditions involving high-order derivatives have been carried out, see for example \cite{EM18} in which a problem of two coupled elliptic equations on a bulk domain and its boundary surface is considered. In this work, both the a priori and a posteriori error estimates for the finite element discretization of the ITP with GTBCs are developed. In particular, residual-type local error estimators are designed to couple the contributions from the bulk equations, transmission boundary conditions and additional tangential flux. Both the upper and local lower bounds of the estimators, which can be used to indicate the reliability and local efficiency of the adaptive FEM, are derived. Numerical experiments for L- and U-shaped domains verify our theoretical results and demonstrate the quasi-optimal convergence of the proposed adaptive FEM. Moreover, considering the significance of the interior transmission eigenvalue problem in inverse scattering, a direct application of the proposed adaptive method to the numerical evaluation of interior transmission eigenvalues is also tested to show the efficiency of the scheme.

The outline of this paper is as follows. Section~\ref{sec:2} introduces the mathematical model of the considered interior transmission problem with GTBCs and establishes the well-posedness of the corresponding variational problem in a nonstandard Sobolev space using the $\mathbb{T}$-coercivity methodology. The finite element approximation together with rigorous a priori and a posteriori error estimates are presented in Section~\ref{sec:3}. Numerical experiments for smooth or piecewise smooth boundaries are shown in Section~\ref{sec:4} to illustrate the accuracy and effectiveness of the proposed adaptive FEM as well as its capability for the computing of interior transmission eigenvalues. Section~\ref{sec:5} finally provides a brief conclusion and future work.

\section{Model problem and well-posedness}
\label{sec:2}
\subsection{Problem description}
\label{sec:2.1}

Let $\Omega\subset\mathbb R^2$ be a nonempty, bounded and connected domain with a piecewise smooth boundary $\Gamma=\bigcup_{j=1}^{N_\Gamma}\Gamma_j$ and each $\Gamma_j$ is smooth. Denote by $\bm \nu$ the unit outer normal to $\Gamma$. This work is interested in the ITP of Helmholtz equations with generalized transmission boundary conditions, which is, given $f\in H^{1/2}(\Gamma)$ and $g\in H^{-1}(\Gamma)$, to determine the wave fields $w$ and $v$ satisfying
\begin{subequations} 
	\label{eq:Strong Model}
	\begin{align}
		&&\Delta w+k^2n(x)w=0 \quad& \text{in}~\Omega, \label{{eq:a}}\\    
		&&\Delta v+k^2v=0 \quad& \text{in}~\Omega,\\
		&&w-v=f \quad& \text{on}~\Gamma,\\    
		&&\rho_1\partial_{\bm{\nu}} w-\rho_2\partial_{\bm{\nu}} v -\mu \Delta_\Gamma w+\lambda w=g \quad& \text{on}~\Gamma,\label{GIBC}\\   &&\nabla_{\Gamma}w|_{\Gamma_j}\cdot\bm{n}_{\partial\Gamma_j}+\nabla_{\Gamma}w|_{\Gamma_k}\cdot\bm{n}_{\partial\Gamma_k}=0	\quad &\text{on}~\partial\Gamma_j\cap\partial\Gamma_k,\label{{eq:d}}
	\end{align} 
\end{subequations}
where $k>0$ is the wave number, $n\in L^\infty(\Omega)$ is the refractive index and $\rho_1,\rho_2>0$. This kind of ITP plays an important role in the study of inverse scattering problems and the considered model relates to the asymptotic model of the acoustic scattering by obstacles with thin layers and coated structures, see for example~\cite{GHK26}. Then the precise constitutive interpretations of $\lambda$ and $\mu$ depend on the underlying wave model and the corresponding asymptotic approximation.
Following~\cite{BCH11}, we assume that $\Re \mu>0$ and $\Im \mu\leq0,\Im \lambda \leq 0$. The signs of $\Im(\lambda)$ and $\Im(\mu)$ are governed by physics to represent absorption, and the assumption on $\Re(\mu)$ is technical to ensure coercivity. In addition, $\Delta_{\Gamma}$ denotes the Laplace-Beltrami operator which is defined by
\begin{equation*}
	\Delta_{\Gamma} u=\nabla_{\Gamma} \cdot \nabla_{\Gamma} u,
\end{equation*}
where $\nabla_{\Gamma}$ is the surface gradient given by
\begin{equation*}
	\nabla_{\Gamma} u=\nabla u-(\nabla u \cdot \bm{\nu}) \bm{\nu}.
\end{equation*} 

It should be pointed out here that of particular interest, values of $k>0$ for which the homogeneous interior transmission problem with $f=g=0$ has a nontrivial solution are called the generalized interior transmission eigenvalues. Analogous to \cite{GHK26,CGH10}, it can be shown that the transmission eigenvalues of the problem (\ref{eq:Strong Model}) form a discrete set accumulating at infinity by assuming additionally $\rho_1>\rho_2,n>1,\Re{\lambda}>0$. Calculation of these particular interior transmission eigenvalues will also be examined in the numerical part. 

\begin{remark}\label{R2.1}
	If the closed surface $\Gamma$ is smooth, the Laplace-Beltrami operator satisfies the integration-by-part formulation on $\Gamma$ as
	\begin{align}\label{eq:integration by parts}
		\int_{\Gamma} \Delta_{\Gamma} u\overline{v} \textup{d} s
		=-\int_{\Gamma} \nabla_{\Gamma} u \cdot \nabla_{\Gamma} \overline{v} \textup{d} s.
	\end{align}
	But for the considered piecewise smooth case, the integration-by-part formulation on each smooth patch $\Gamma_j$ shall be modified as
	\begin{equation*}
		\int_{\Gamma_j} \Delta_{\Gamma} u\overline{v} \textup{d} s
		=-\int_{\Gamma_j} \nabla_{\Gamma} u \cdot \nabla_{\Gamma} \overline{v} \textup{d} s+ \sum_{\bm{x} \in \partial \Gamma_j} \nabla_{\Gamma} u(\bm{x}) \cdot \bm{n}_{\partial \Gamma_j} \overline{v}(\bm{x}),
	\end{equation*}
	where $\bm{n}_{\partial \Gamma_j}$ denotes the outer unit normal of $\partial \Gamma_j$; see also \cite{KCDQ15}. That is why we need an additional compatibility condition (\ref{eq:Strong Model}e) here for the purpose of studying the well-posedness of the problem (\ref{eq:Strong Model}) via a variational approach.
\end{remark}

\subsection{Variational Problem and well-posedness analysis}
\label{sec:2.2}

It is clear that in the well-posedness analysis of (\ref{eq:Strong Model}), a difficulty is encountered due to the transmission condition, which involves higher-order derivatives on $\Gamma$ induced by the Laplace-Beltrami operator. Consequently, the $H^{1/2}(\Gamma)$-regularity provided by the standard trace theorem for $H^1$ Sobolev space falls short of defining a proper weak formulation of problem (\ref{eq:Strong Model}). For the piecewise smooth boundary $\Gamma$, we define
\begin{equation*}
	H^m(\Gamma):=\left\{w\in L^2(\Gamma):w|_{\Gamma_j}\in H^m(\Gamma_j),\; [w]=0\ \text{on } \partial\Gamma_j\cap\partial\Gamma_k\;\;\text{for all}\;\; j,k=1,\cdots,N_\Gamma \right\},\quad m > 0,
\end{equation*}
equipped with the norm
\begin{equation*}
	\|w\|_{H^m(\Gamma)}^2:=\sum_{j=1}^{N_\Gamma}\|w\|_{H^m(\Gamma_j)}^2.
\end{equation*}
Here, $[w]$ denotes the jump of the function $w$ across the corners $\partial\Gamma_j\cap\partial\Gamma_k$. Then we define a nonstandard Sobolev space
\begin{equation*}
	H^m(\Omega;\Gamma):=\left\{w\in H^m(\Omega):w|_{\Gamma}\in H^m(\Gamma)\right\},\quad m > 0,
\end{equation*}
For convenience of the subsequent analysis, define the product Sobolev spaces
\begin{align*}
	&\mathcal{H}^m(\Omega;\Gamma):=H^m(\Omega;\Gamma)\times H^m(\Omega),\quad m > 0,\\
	&\mathcal{H}_{\Gamma}^m(\Omega;\Gamma):=\left\{(w, v) \in \mathcal{H}^m(\Omega;\Gamma): w=v ~ \text{on} ~\Gamma\right\},\quad m > 0,
\end{align*}
both equipped with the norm
\begin{equation*}
	\|\bm{U}\|_{\mathcal{H}^m(\Omega;\Gamma)}=
	\left(\|w\|_{H^m(\Omega)}^2+\|w\|_{H^m(\Gamma)}^2+\|v\|_{H^m(\Omega)}^2\right)^{1 / 2}, \quad \forall \bm{U}=(w, v) \in \mathcal{H}^m(\Omega;\Gamma).
\end{equation*}
In addition, we define
\begin{align*}
	L^2(\Omega;\Gamma):=L^2(\Omega)\times L^2(\Gamma) \times L^2(\Omega),
\end{align*}
equipped with the norm
\begin{align*}
	\|\bm{U}\|_{L^2(\Omega;\Gamma)}^2:=
	\|w\|_{L^2(\Omega)}^2+\|w\|_{L^2(\Gamma)}^2+\|v\|_{L^2(\Omega)}^2.
\end{align*}

The variational formulation of the transmission problem (\ref{eq:Strong Model}) reads: given $f\in H^{1/2}(\Gamma)$ and $g\in H^{-1}(\Gamma)$, to find  $\bm{U}=(w,v) \in \mathcal{H}^1(\Omega;\Gamma)$ such that $w - v = f$ on $\Gamma$ and
\begin{equation}\label{eq:variational problem}
	A(\bm{U}, \Psi)=\ell(\Psi),\quad ~\forall \Psi=(\psi;\phi)\in \mathcal{H}_{\Gamma}^1(\Omega;\Gamma),
\end{equation} 
where 
\begin{align}\label{abc}
	A(\bm{U}, \Psi)= &b(\bm{U}, \Psi)-c(\bm{U}, \Psi),
\end{align}
with
\begin{align*}
	b(\bm{U}, \Psi) & =\rho_1\int_\Omega(\nabla w\cdot \nabla \overline{\psi}+w\overline{\psi})\text{d} \bm{x}-\rho_2\int_\Omega(\nabla v\cdot\nabla \overline{\phi}+v\overline{\phi})\text{d} \bm{x}+\int_\Gamma(\mu \nabla_\Gamma w\cdot\nabla_\Gamma\overline{\psi}+\delta w\overline{\psi})\text{d} s,\\    
	c(\bm{U}, \Psi) & =\rho_1\int_\Omega(k^2nw\overline{\psi}+w\overline{\psi})\text{d} \bm{x}-\rho_2\int_\Omega(k^2v\overline{\phi}+v\overline{\phi})\text{d} \bm{x}+\int_\Gamma(\delta w\overline{\psi}-\lambda w\overline{\psi})\text{d} s, 
\end{align*}
and where $\ell:H^{-1}(\Gamma) \rightarrow \mathbb{C}$ is a linear functional 
\begin{equation} \label{eq:function_ell}
	\ell(\Psi)  =\int_{\Gamma} g \overline{\psi} \text{d} s.
\end{equation}
Here, the parameter $\delta>0$ is artificial and its value will be specified later. In particular, the Sobolev embedding theorem implies that $c(\cdot, \cdot)$ is a compact perturbation of $b(\cdot, \cdot)$.
However, $b(\cdot, \cdot)$ is not coercive here, and thus the Fredholm framework cannot be used directly. To overcome this difficulty and extend the Fredholm alternative theorem to the
variational problem (\ref{eq:variational problem}), the so-called $\mathbb{T}$-coercivity approach will be used in a way similar to \cite{WC13}.

\begin{definition}\label{def：T-cocecivity}
	Let $V$ be a Hilbert space and let $\mathbb{T}:\bm{V} \rightarrow \bm{V}$ be a continuous linear operator in $\bm{V}$. A sesquilinear form $b$ is said to be $\mathbb{T}$-coercive in $\bm{V} \times \bm{V}$, if there exists some $\delta_T > 0$ such that
	\begin{equation*}
		{\left|b\left(\Psi,\mathbb{T}\Psi \right)\right|}\geq\delta_T \|\Psi \|^2_{\bm{V}},\quad \forall \Psi \in \bm{V}.
	\end{equation*}
\end{definition} 

{Let $(\cdot,\cdot)_{1,\Omega}$} denote the inner product of the space $H^1(\Omega)$. We first state the following lifting property. {For brevity, the notations $a \lesssim b$ and $a \gtrsim b $ are introduced to imply that $a \leq Cb$ and $a \geq Cb$, respectively.}

\begin{lemma}\label{lem:trace lifting}
	For every $q\in H^{1/2}(\Gamma)$, there exists a unique $\mathcal{L}q\in H^{1}(\Omega)$ such that
	\begin{equation*}
		\mathcal{L}q|_\Gamma=q,\qquad(\mathcal{L}q,z)_{1,\Omega}=0, \quad \forall z\in H^1_0(\Omega),
	\end{equation*}
	and the lifting operator  $\mathcal{L}:H^{1/2}(\Gamma)\rightarrow H^1(\Omega)$ is bounded, i.e.,
	\begin{equation}\label{ieq:lifting_estimates}
		\|\mathcal{L}q\|_{H^1(\Omega)}\lesssim \|q\|_{H^{1/2}(\Gamma)}.
	\end{equation}
\end{lemma}
\begin{proof}
	This result follows directly from the trace theorem and the well-posedness theory for the Dirichlet problem of modified Helmholtz equation $\Delta v- v=0$ in $\Omega$, See,e.g. \cite[Chapter~3,Theorem~3.37]{M00}. Hence the proof is omitted here.
\end{proof}

	

Since $H^1(\Gamma)$ is compactly embedded in $H^{1/2}(\Gamma)$, which is continuously embedded in $L^2(\Gamma)$. The estimate (\ref{ieq:lifting_estimates}) together with Ehrling's inequality (also known as the J.-L. Lions lemma\cite[p.173]{B11}) further imply that, for any $q\in H^1(\Gamma)$ and any $\alpha>0$, there exist a unique $\mathcal{L}q\in H^{1}(\Omega)$ and a constant $C_\alpha>0$ such that
\begin{equation}\label{ieq:lifting estimates with ehrling}
	\|\mathcal{L}q\|^2_{H^1(\Omega)}\leq \alpha\|q\|^2_{H^{1}(\Gamma)}+C_\alpha\|q\|^2_{L^{2}(\Gamma)}.
\end{equation}

Now we define $\mathbb{T}:\mathcal{H}_\Gamma^1(\Omega;\Gamma)\rightarrow\mathcal{H}_\Gamma^1(\Omega;\Gamma)$ as
\begin{equation*}
	\mathbb{T}\Psi=(\psi,2\mathcal{L}q-\phi),\quad \forall \Psi=(\psi;\phi)\in \mathcal{H}_\Gamma^1(\Omega;\Gamma),
\end{equation*}
where $\mathcal{L}q|_\Gamma =q=\phi|_\Gamma=\psi|_\Gamma$.
It is easy to check that $\mathbb{T}$ is bounded and bijective since $\mathbb{T}^2=I$.

\begin{lemma}\label{lem:btcocecive}
	The sesquilinear form $b(\cdot, \cdot)$ is $\mathbb{T}$-coercive in $\mathcal{H}^1_\Gamma(\Omega;\Gamma)\times\mathcal{H}^1_\Gamma(\Omega;\Gamma)$.
\end{lemma}
\begin{proof}
	For any $\Psi=(\psi;\phi) \in \mathcal{H}^1_\Gamma(\Omega;\Gamma)$ with $q=\phi|_\Gamma=\psi|_\Gamma$,
	\begin{align*}
		\Re{b(\Psi, \mathbb{T}\Psi)} =&\rho_1\int_\Omega(\nabla \psi\cdot\nabla \overline{\psi}+\psi\overline{\psi})\text{d} \bm{x}+\rho_2\int_\Omega(\nabla \phi\cdot\nabla \overline{\phi}+\phi\overline{\phi})\text{d} \bm{x}+\int_\Gamma(\Re{\mu} \nabla_\Gamma \psi\cdot\nabla_\Gamma\overline{\psi}+\delta \psi\overline{\psi})\text{d} s\\&-2\rho_2\int_\Omega(\nabla \phi\cdot\nabla \overline{\mathcal{L}q}+\phi\overline{\mathcal{L}q})\text{d} \bm{x}.
	\end{align*}
	Applying Young’s inequality and using (\ref{ieq:lifting estimates with ehrling}), we obtain
	\begin{align*}
		\Re{b(\Psi, \mathbb{T}\Psi)} &\geq \rho_1\|\psi\|^2_{H^1(\Omega)}+\rho_2\|\phi\|^2_{H^1(\Omega)}+\Re{\mu} \|\psi\|^2_{H^1(\Gamma)}+\left(\delta-\Re{\mu}\right)\|\psi\|^2_{L^2(\Gamma)}\\&\quad -2\rho_2\left(\beta\|\phi\|^2_{H^1(\Omega)}+C_\beta\|\mathcal{L}q\|^2_{H^1(\Omega)}\right)\\
		&\geq \rho_1\|\psi\|^2_{H^1(\Omega)}+\rho_2(1-2\beta)\|\phi\|^2_{H^1(\Omega)}+\Re{\mu} \|\psi\|^2_{H^1(\Gamma)}+(\delta-\Re{\mu})\|\psi\|^2_{L^2(\Gamma)}\\&\quad -2\rho_2C_\beta\left(\alpha\|q\|^2_{H^1(\Gamma)}+C_\alpha\|q\|^2_{L^2(\Gamma)}\right),
	\end{align*}
	which means that
	\begin{equation*}
		\Re{b(\Psi, \mathbb{T}\Psi)} \geq \rho_1\|\psi\|^2_{H^1(\Omega)}+\rho_2(1-2\beta)\|\phi\|^2_{H^1(\Omega)}+\left(\Re{\mu}-2\rho_2\alpha C_\beta \right)\|\psi\|^2_{H^1(\Gamma)}+\left(\delta-\Re{\mu}-2\rho_2 C_\alpha C_\beta \right)\|\psi\|^2_{L^2(\Gamma)}.
	\end{equation*} 
	Then the proof is completed by choosing {$0<\beta<1/2$, $0<\alpha<\Re\mu/(2\rho_2C_\beta)$}, and $\delta>\Re{\mu}+2\rho_2 C_\alpha C_\beta >0$ in sequence, and letting $\delta_T=\text{min}\{ \rho_1,\rho_2(1-2\beta),\Re{\mu}-2\rho_2\alpha C_\beta\}>0$.
\end{proof}

\begin{theorem}\label{Th-vpunique}
	Assume that $k$ is not a transmission eigenvalue. Then for every $f\in H^{1/2}(\Gamma)$, $g\in H^{-1}(\Gamma)$, there exists a unique solution $\bm{U}=(w,v) \in \mathcal{H}^1(\Omega;\Gamma)$ satisfying the variational problem (\ref{eq:variational problem}). In addition, the inf-sup condition
	\begin{equation}\label{ieq:vpinf-sup}
		\sup_{0\ne\Psi\in{\mathcal{H}^1_\Gamma(\Omega;\Gamma)}}
		\frac{|A(\bm{U}_0,\Psi)|}{\|\Psi\|_{\mathcal{H}^1(\Omega;\Gamma)}}\gtrsim \|\bm{U}_0\|_{\mathcal{H}^1(\Omega;\Gamma)},
		\quad \forall \bm{U}_0 \in \mathcal{H}^1_\Gamma(\Omega;\Gamma).
	\end{equation}
\end{theorem}
\begin{proof}
	Since $\mathbb{T}$ is bijective, the variational problem (\ref{eq:variational problem}) is equivalent to the problem of finding $\bm{U}=(w,v) \in \mathcal{H}^1(\Omega;\Gamma)$ such that $w-v=f$ on $\Gamma$, and
	\begin{equation}\label{eq:uniuque}
		b(\bm{U}, \mathbb{T}\Psi)-c(\bm{U}, \mathbb{T}\Psi)=\ell(\mathbb{T}\Psi),\quad ~\forall \Psi=(\psi;\phi)\in \mathcal{H}_{\Gamma}^1(\Omega;\Gamma).
	\end{equation}
	Take $\Theta_f =(0,-\mathcal{L}f)$, equation (\ref{eq:uniuque}) can be rewritten as: find $\bm{U}=\bm{U}_0+\Theta_f,\bm U_0\in \mathcal{H}^1_\Gamma(\Omega;\Gamma)$ such that
	\begin{equation}\label{eq:uniuque1}
		A(\bm{U}_0, \mathbb{T}\Psi)=b(\bm{U}_0, \mathbb{T}\Psi)-c(\bm{U}_0, \mathbb{T}\Psi)=\ell^*(\mathbb{T}\Psi),\quad ~\forall \Psi=(\psi;\phi)\in \mathcal{H}_{\Gamma}^1(\Omega;\Gamma),
	\end{equation}
	where
	\begin{equation}\label{eq:T dual product}
		\ell^*(\mathbb{T}\Psi)= \ell(\mathbb{T}\Psi)-A(\Theta_f, \mathbb{T}\Psi).
	\end{equation} 
	It can be shown that $\ell^*$ is bounded as
	\begin{align*}
		|\langle \ell^*(\Psi) \rangle|&\leq \|g\|_{H^{-1}(\Gamma) }\|\psi\|_{H^{1}(\Gamma) }+\rho_2(1+k^2)\|\mathcal{L}f\|_{H^{1}(\Omega)}\|\phi\|_{H^{1}(\Omega)}\\
		&\leq \left(\|g\|_{H^{-1}(\Gamma)}+\rho_2(1+k^2)C_0\|f \|_{H^{1/2}(\Gamma)}\right)\|\Psi\|_{\mathcal{H}^1(\Omega;\Gamma)}.
	\end{align*}
	\Cref{lem:btcocecive} and the compactness of $c(\cdot,\mathbb{T}\cdot)$ imply that the Riesz operator induced by $A(\cdot,\mathbb{T}\cdot)$ is Fredholm of index zero. Then the assumption that $k$ is not a transmission eigenvalue yields the well-posedness of the variational problem (\ref{eq:uniuque1}) by Fredholm alternative. The proof is complete.
\end{proof}

\section{Finite element analysis}
\label{sec:3}

In this section, we consider the finite element discretization of the variational problem (\ref{eq:variational problem}) and derive both the a priori and a posteriori error estimates. For simplicity, the discussion of the approximation error between the boundaries $\Gamma$ and the finite element mesh, which is beyond the main scopes of this work, will be ignored. For readers interested in the analysis of iso-parametric element method for problems with the generalized impedance boundary condition, we refer to \cite{CGP24}.

\subsection{Finite element discretization}
\label{sec:3.1}
Let 
$\mathcal{M}_{h}$ be the regular triangulation  of $\Omega$, where $h$ denotes the mesh size. Assume that $X_{h} \subset H^1(\Omega)$ is the standard conforming finite element space 
\begin{align*}
	X_{h}^p&:=\{ w_h \in C(\overline{\Omega}): w_h|_T \in P_p(T) ~\text{for all} ~T \in \mathcal{M}_{h} \},	\\
	X_{h}&:=X_{h}^1,
\end{align*} 
where $P_p(T)$ denotes the set of all polynomials of degree no more than $p$.   In addition, the boundary mesh inherited from $\mathcal{M
}_h$ is set as
\begin{equation*}
	\mathcal{S}_h	=\{ S \in \Gamma; ~\textup{there exists some}~ T \in \mathcal{M}_h ~\textup{such that}~ S=\overline{T}\cap\Gamma  \}.
\end{equation*}

For convenience of the subsequent error analysis, introduce the product finite element spaces 
\begin{align*}
	&\mathcal{H}_h := X_{h} \times X_{h} \subset \mathcal{H}^1(\Omega;\Gamma),\\
	&\mathcal{H}_{\Gamma,h}=\left\{\left(w_h, v_h\right) \in \mathcal{H}_h: w_h=v_h~ \text{on}~ \Gamma\right\}.
\end{align*}
Let $X_{\Gamma,h}$ denote the restriction of $X_h$ to $\Gamma$ which consists of piecewise linear functions on $\Gamma$. Define the orthogonal projection $Q^h: L^2(\Gamma) \rightarrow X_{\Gamma,h}$ as follows
$$
\left\langle Q^h f, \chi\right\rangle=\langle f, \chi\rangle, \quad \forall \chi \in X_{\Gamma,h}, 
$$
where $\langle\cdot, \cdot\rangle$ denotes the inner product on $L^2(\Gamma)$. It follows that for any $f \in H^r(\Gamma)$,
$$
\begin{aligned}
	\left\|f-Q^h f\right\|_{H^{-1 / 2}(\Gamma)} & \lesssim h^{r+1 / 2}\|f\|_{H^r(\Gamma)}, \quad 0<r<2,\\
	\left\|f-Q^h f\right\|_{L^2(\Gamma)} & \lesssim h^r\|f\|_{H^r(\Gamma)}, \quad 0<r<2, \\
	\left\|f-Q^h f\right\|_{H^{1 / 2}(\Gamma)} & \lesssim h^{r-1 / 2}\|f\|_{H^r(\Gamma)}, \quad 0<r<2,
\end{aligned}
$$
which further implies that
\begin{equation}\label{ieq:projection}
	\left\|f-Q^h f\right\|_{H^{-s}(\Gamma)} \lesssim h^{r+s}\|f\|_{H^r(\Gamma)}, \quad 0<r<2,~-1/2 \leq s \leq 1 / 2,
\end{equation}
by interpolation. We also define the discrete lifting operator $\mathcal{L}^h:X_{\Gamma,h}\rightarrow X_h,q^h\rightarrow \mathcal{L}^hq^h$ such that 
\begin{equation*}        \mathcal{L}^hq^h|_\Gamma=q^h,\qquad(\mathcal{L}^hq^h,z^h)_{1,\Omega}=0, \quad \forall z^h\in X_h\cap H^1_0(\Omega).
\end{equation*}
Then
\begin{equation}\label{ieq:discrete_lifting_estimates}
	\|\mathcal{L}^hq^h\|_{H^1(\Omega)}\lesssim\|q^h\|_{H^{1/2}(\Gamma)},
\end{equation}
The existence and uniform stability of such discrete lifting operator have been proved in \cite{DS03} by using the boundary-preserving Scott–Zhang projection. 

Letting $f^h:=Q^hf$, the finite element approximation to the variational problem \eqref{eq:variational problem} reads as follows: find $\bm{U}^h=(w^h,v^h) \in \mathcal{H}_h$ such that $w^h-v^h=f^h$ on $\Gamma$, and 
\begin{equation}\label{eq:discrete problem}
	A(\bm{U}^h, \Psi^h)=\ell(\Psi^h),\quad \forall \Psi^h=(\psi^h;\phi^h)\in \mathcal{H}_{\Gamma ,h}.
\end{equation}
Combining both the $\mathbb{T}$-coercivity of $b$ on $\mathcal{H}_{\Gamma,h}\times \mathcal{H}_{\Gamma,h}$ and a general argument of Schatz \cite{S74} subsequently validates that, for sufficiently small $h>0$, there holds the inf-sup condition:
\begin{equation}\label{ieq:dpinf-sup}
	\sup_{0\ne\Psi^h\in{\mathcal{H}_{\Gamma,h}}}
	\frac{|A(\bm{U},\Psi^h)|}{\|\Psi^h\|_{\mathcal{H}^1(\Omega;\Gamma)}}\gtrsim \|\bm{U}\|_{\mathcal{H}^1(\Omega;\Gamma)},
	\quad \forall \bm{U} \in \mathcal{H}_{\Gamma,h}.
\end{equation}
For brevity, the detailed derivation is omitted here. One can also refer to \cite{BBCZ10} for an analogous proof. 


\subsection{A priori error estimates}
\label{sec:3.2}

We first establish the a priori error estimates in the energy norm $\mathcal{H}^1(\Omega;\Gamma)$. The following Cea's lemma is required.

\begin{lemma}\label{th:inf H1 priori estimates}
	{Let $\bm{U}=(w,v) \in \mathcal{H}^1(\Omega;\Gamma)$ and $\bm{U}^h=(w^h,v^h) \in \mathcal{H}_h$ be the solutions of \eqref{eq:variational problem} and \eqref{eq:discrete problem}, respectively. Then there exists a constant $h_0>0$ such that for any $0<h \leq h_0$, the following estimate holds}
	\begin{align*}
		\left\|w-w^h\right\|_{H^1(\Omega)}+\left\|v-v^h\right\|_{H^1(\Omega)} +\left\|w-w^h\right\|_{H^1(\Gamma)}\lesssim 
		& \inf_{\hat{w}^h\in X^h}\left(\left\|w-\hat{w}^h\right\|_{H^1(\Omega)}+\left\|w-\hat{w}^h\right\|_{H^1(\Gamma)}\right)\\
		& +\inf_{\hat{v}^h\in X^h}\left\|v-\hat{v}^h\right\|_{H^1(\Omega)} 
		+\left\|f-f^h\right\|_{H^{1/2}(\Gamma)}
		.
	\end{align*}
\end{lemma}
\begin{proof}
	For any $ \widetilde {\bm{U}}^h=( \widetilde{w}^h,\widetilde{v}^h) \in \mathcal{H}_h$ such that $ \widetilde{w}^h-\widetilde{v}^h=f^h$, it is easy to verify that $\widetilde{\bm{U}}^h-\bm{U}^h \in \mathcal{H}_{\Gamma,h}$. Then the discrete inf-sup
	condition \eqref{ieq:dpinf-sup} and the Galerkin orthogonality from \eqref{eq:variational problem} and \eqref{eq:discrete problem} yield that
	\begin{align*}
		\left\|\widetilde{\bm{U}}^h-\bm{U}^h\right\|_{\mathcal{H}^1(\Omega; \Gamma)}
		&\lesssim \sup_{0\ne\Psi^h\in{\mathcal{H}_{\Gamma,h}}}
		\frac{|A(\widetilde{\bm{U}}^h-\bm{U}^h,\Psi^h)|}{\|\Psi^h\|_{\mathcal{H}^1(\Omega; \Gamma)}}\\
		&= \sup_{0\ne\Psi^h\in{\mathcal{H}_{\Gamma,h}}}
		\frac{|A(\bm{U}-\widetilde{\bm{U}}^h,\Psi^h)|}{\|\Psi^h\|_{\mathcal{H}^1(\Omega; \Gamma)}}
		\lesssim  \left\|\bm{U}-\widetilde{\bm{U}}^h\right\|_{\mathcal{H}^1(\Omega; \Gamma)}.
	\end{align*}
	Then, by the triangle inequality and passing to the infimum, it follows 
	\begin{equation}\label{ieq:best-fit optimal approximation}
		\left\|\bm{U}-\bm{U}^h\right\|_{\mathcal{H}^1(\Omega; \Gamma)}\lesssim \inf_{\substack{0\ne\widetilde{\bm{U}}^h\in{\mathcal{H}_{h}} \\\widetilde{w}^h-\widetilde{v}^h=f^h}}\left\|\bm{U}-\widetilde{\bm{U}}^h\right\|_{\mathcal{H}^1(\Omega; \Gamma)}.
	\end{equation}
	Fix an arbitrary $\hat{w}^h \in X_h$, we define $p^h:=\hat{w}^h|_\Gamma-f^h$. For any $\hat{v}^h_0\in X_h$, set $\hat{v}^h:=\hat{v}^h_0+\mathcal{L}^h(p^h-\hat{v}^h_0|_\Gamma)$.
	It is easy to know $\hat{v}^h|_\Gamma=p^h$ and therefore $\hat{w}^h|_\Gamma-\hat{v}^h|_\Gamma=f^h$.
	By using the triangle inequality, \eqref{ieq:discrete_lifting_estimates},  and the trace theorem yield that
	\begin{align*}
		\left\|v-\hat{v}^h\right\|_{H^1(\Omega)}
		&\leq\left\|v-\hat{v}^h_0\right\|_{H^1(\Omega)}+\left\|\hat{v}^h-\hat{v}^h_0\right\|_{H^1(\Omega)}\\
		&\lesssim   \left\|v-\hat{v}^h_0\right\|_{H^1(\Omega)}+ \left\|p^h-\hat{v}^h_0\right\|_{H^{1/2}(\Gamma)}\\
		&\lesssim \left\|v-\hat{v}^h_0\right\|_{H^1(\Omega)}+\left\|p^h-v\right\|_{H^{1/2}(\Gamma)}.
	\end{align*}
	Since $v|_\Gamma=w|_\Gamma-f$, hence $p^h-v|_\Gamma=(\hat{w}^h|_\Gamma-w|_\Gamma)+(f-f^h)$.
	Using the continuous embedding $H^1(\Gamma) \hookrightarrow H^{1/2}(\Gamma)$, the above inequality further implies 
	\begin{align}\label{ieq:v best-fit optimal approximation}
		\left\|v-\hat{v}^h\right\|_{H^1(\Omega)}
		&\lesssim  \left\|v-\hat{v}^h_0\right\|_{H^1(\Omega)}+\left\|w-\hat{w}^h\right\|_{H^{1/2}(\Gamma)}+\|f-f^h\|_{H^{1/2}(\Gamma)} \nonumber \\    
		&\lesssim  \left\|v-\hat{v}^h_0\right\|_{H^1(\Omega)}+\left\|w-\hat{w}^h\right\|_{H^1(\Gamma)}+\|f-f^h\|_{H^{1/2}(\Gamma)}.
	\end{align}
	Let $\widetilde{\bm{U}}^h=(\hat{w}^h,\hat{v}^h)$. Then substituting \eqref{ieq:v best-fit optimal approximation} into \eqref{ieq:best-fit optimal approximation}, and taking the infimum over $\hat{w}^h,\hat{v}^h_0 \in X_h$ proves the result.   
\end{proof}

Denote by $\Pi_\Omega^h:H^1(\Omega)\rightarrow X_h,\Pi_\Gamma^h:H^1(\Gamma)\rightarrow X_{\Gamma,h}$ the Scott–Zhang quasi-interpolation operators associated with the bulk mesh and the induced boundary mesh, respectively. Their standard approximation properties, including fractional-order estimates, can be found in \cite{SZ90}. In particular, for $w \in H^{1+r_w}(\Omega),w|_\Gamma \in H^{1+s_\Gamma}(\Gamma),0<r_w,s_\Gamma\leq1$, the following estimates hold:
\begin{equation}\label{ieq:bulk quasi-interpolation}
	\|w-\Pi_\Omega^h w\|_{H^1(\Omega)}\lesssim h^{r_w}\|w\|_{H^{1+r_w}(\Omega)},
\end{equation}
and
\begin{equation}\label{ieq:surface quasi-interpolation}
	\|w-\Pi_\Gamma^h (w|_\Gamma)\|_{H^s(\Gamma)}\lesssim h^{1+s_\Gamma-s}\|w\|_{H^{1+s_\Gamma}(\Gamma)}, \quad 1/2\leq s \leq1.
\end{equation}
The fractional-order estimates follow from the standard approximation properties and Sobolev interpolation inequalities.

\begin{lemma}\label{lem:Bulk–surface compatible quasi-interpolation}
	Assume that $w \in H^{1+r_w}(\Omega)$ and $w|_\Gamma \in H^{1+s_\Gamma}(\Gamma)$ with $0<r_w,s_\Gamma\leq1$.
	Define
	\begin{equation}\label{eq:bs quasi-interpolation}
		\widetilde{\Pi}_\Omega^h w  :=\Pi_\Omega^h w+\mathcal{L}^h(\Pi_\Gamma^h(w|_\Gamma)-(\Pi_\Omega^h w)|_\Gamma).
	\end{equation}
	Then 
	\begin{equation}\label{ieq:bulk estimates bs quasi interplation}
		\|w-\widetilde{\Pi}_\Omega^h w\|_{H^1(\Omega)}\lesssim  h^{r_w}\|w\|_{H^{1+r_w}(\Omega)}+h^{1/2+s_\Gamma}\|w\|_{H^{1+s_\Gamma}(\Gamma)},
	\end{equation}
	and
	\begin{equation}\label{ieq:surface estimates bs quasi interplation}
		\|w-(\widetilde{\Pi}_\Omega^h w)|_\Gamma\|_{H^1(\Gamma)}\lesssim h^{s_\Gamma}\|w\|_{H^{1+s_\Gamma}(\Gamma)}.
	\end{equation}

\end{lemma}
\begin{proof}
	Since both $\Pi_\Gamma^h(w|_\Gamma)$ and $(\Pi_\Omega^h w)|_\Gamma$  belong to $X_{\Gamma,h}$, the operator in \eqref{eq:bs quasi-interpolation} is well defined. From the lifting property of $\mathcal{L}^h$, 
	\begin{equation}\label{eq:wbs quasi-interpolation}
		(\widetilde{\Pi}_\Omega^h  w)|_\Gamma=(\Pi_\Omega^h w)|_\Gamma+\left(\mathcal{L}^h(\Pi_\Gamma^h(w|_\Gamma)-(\Pi_\Omega^h w)|_\Gamma)\right)|_\Gamma=\Pi_\Gamma^h(w|_\Gamma).
	\end{equation}
	Using the boundedness of $\mathcal{L}^h$, the triangle inequality, and the continuity of the trace operator yields
	\begin{align*}
		\|w-\widetilde{\Pi}_\Omega^h w\|_{H^1(\Omega)}
		&\leq\|w-\Pi_\Omega ^h w\|_{H^1(\Omega)}+\|\mathcal{L}^h(\Pi_\Gamma^h(w|_\Gamma)-(\Pi_\Omega^h w)|_\Gamma)\|_{H^{1}(\Omega)}\\
		&\lesssim \|w-\Pi_\Omega ^h w\|_{H^1(\Omega)}+\|\Pi_\Gamma^h(w|_\Gamma)-(\Pi_\Omega^h w)|_\Gamma\|_{H^{1/2}(\Gamma)}\\
		&\leq\|w-\Pi_\Omega ^h w\|_{H^1(\Omega)}+\|w-(\Pi_\Omega^h w)|_\Gamma\|_{H^{1/2}(\Gamma)}+\|w-\Pi_\Gamma^h(w|_\Gamma)\|_{H^{1/2}(\Gamma)}\\
		&\lesssim \|w-\Pi_\Omega ^h w\|_{H^1(\Omega)}+\|w-\Pi_\Gamma^h(w|_\Gamma)\|_{H^{1/2}(\Gamma)}.
	\end{align*}
	Applying \eqref{ieq:bulk quasi-interpolation} and \eqref{ieq:surface quasi-interpolation} with $s=\frac{1}{2}$ to the above inequality gives \eqref{ieq:bulk estimates bs quasi interplation}.
	Finally, taking $s=1$ and combining \eqref{ieq:surface quasi-interpolation} with \eqref{eq:wbs quasi-interpolation} shows \eqref{ieq:surface estimates bs quasi interplation}. The proof is complete. 
\end{proof}

Now we take $\hat{w}^h=\widetilde{\Pi}_\Omega^h w,\hat{v}^h=\Pi_\Omega^hv$, we obtain the following a priori error estimates by applying \eqref{ieq:projection}, \eqref{ieq:bulk quasi-interpolation}, \eqref{ieq:bulk estimates bs quasi interplation} and \eqref{ieq:surface estimates bs quasi interplation} to Lemma~\ref{th:inf H1 priori estimates}. 

\begin{theorem}\label{th:H1 error}
	Let $\bm{U}$ and $\bm{U}^h$ be the solutions of \eqref{eq:variational problem} and  \eqref{eq:discrete problem}, respectively. Suppose that $w \in H^{1+r_w}(\Omega),v \in H^{1+r_v}(\Omega),w|_\Gamma \in H^{1+s_\Gamma}(\Gamma)$, and $ f \in H^{1/2+s_f}(\Gamma)$, where $0<r_w,r_v,s_\Gamma, s_f\leq1$. Then there exist a constant $h_0>0$ such that for any $0<h \leq h_0$, the following estimate holds
	\begin{align*}
		\left\|\bm{U}-\bm{U}^h\right\|_{\mathcal{H}^1(\Omega;\Gamma)}\lesssim 
		& h^{r_w} \|w\|_{H^{1+r_w}(\Omega)}+h^{r_v} \|v\|_{H^{1+r_v}(\Omega)}\\
		+& h^{s_\Gamma} \|w\|_{H^{1+s_\Gamma}(\Gamma)}+h^{s_f}\|f\|_{H^{1/2+s_f}(\Gamma)}.
	\end{align*}
\end{theorem}

\subsection{A posteriori error estimates: upper bound}
\label{sec:3.3}

The a priori error estimate shows that the linear finite element method on quasi-uniform triangulations can achieve the optimal convergence rate $\mathcal{O}(h)$ under the regularity assumption $\bm{U}\in \mathcal{H}^2(\Omega;\Gamma)$. For domains with non-smooth boundaries, however, this regularity assumption may fail to hold. In such cases, adaptive finite element techniques may be employed to enhance numerical accuracy and computational efficiency.

To begin with, we define the a posteriori error indicator in terms of several residuals for the problem \eqref{eq:variational problem}, which is the fundamental component of the adaptive finite element method. Inside the domain $\Omega$, the residual on each element $T \in \mathcal{M}_{h}$ is defined by
\begin{align}
	R_{w,T}&= \rho_1 (\Delta w^h + k^2nw^h) \label{eq:Rw},\\
	R_{v,T}&= \rho_2 ( \Delta v^h + k^2v^h) \label{eq:Rv}.
\end{align}
Let $\partial T$ be the set of all edges of $T$, and denote the length of an edge $e \in \partial T$ by $h_e$. For any interior edge $e$ which is the common side of $T_1,~ T_2\in  \mathcal{M}_{h}$, define the jump residual across $e$ as 
\begin{align}
	J_{w,e}&=\rho_1\left(\nabla w^h|_{T_1} \cdot \bm{n}_1+\nabla w^h|_{T_2} \cdot \bm{n}_2 \right),\label{eq:Jw}\\
	J_{v,e}&=\rho_2\left(\nabla v^h|_{T_1} \cdot \bm{n}_1+\nabla v^h|_{T_2} \cdot \bm{n}_2 \right) \label{eq:Jv},
\end{align}
where $\bm{n}_i$ is the unit outward normal vector to the boundary of $T_i,~i=1,2$. Corresponding to the transmission boundary condition (\ref{eq:Strong Model}d), the residual on each boundary element $S\in \mathcal{S}_h$ is defined by
\begin{equation}
	R_{\Gamma,S}= g+\mu\Delta_{\Gamma}w^h-\lambda w^h- \rho_1\partial_{\bm \nu} w^h+\rho_2\partial_{\bm \nu} v^h\label{eq:Rx}.
\end{equation}
And for any $\bm{x}=\partial S_1 \cap\partial S_2$ which is the common endpoint of $S_1,S_2 \in S_h$, define
\begin{align}
	J_{\bm{x}}=\mu\left(\nabla_{\Gamma}w^h \cdot \bm{n}_{\partial S_1}+\nabla_{\Gamma}w^h\cdot \bm{n}_{\partial S_2}\right)\label{eq:Jx},
\end{align}
where $\bm{n}_{\partial S_i}$ is the unit outer normal vector of  $\partial S_i, i=1,2$. This residual relates to the compatibility condition (\ref{eq:Strong Model}e). For $T \in \mathcal{M}_{h}$ and $S \in \mathcal{S}_{h}$, now we introduce the local error estimators as
\begin{equation} \label{eq:local estimator}   \eta_{T}^2=\eta_{w,T}^2+\eta_{v,T}^2+\sum_{S\subset \partial T \cap\Gamma} \eta_{S}^2,
\end{equation} 
where
\begin{align}
	&\eta_{w,T}^2=h_T^2\|  R_{w,T}\|_{L^2(T)}^2+	\left( \frac{1}{2} \sum_{e \in \partial T \setminus \Gamma} h_e\| J_{w,e}\|_{L^2(e)}^2 \right),\label{eq:local estimator w}\\
	&\eta_{v,T}^2=h_T^2\|  R_{v,T}\|_{L^2(T)}^2+	\left( \frac{1}{2} \sum_{e \in \partial T \setminus \Gamma} h_e\| J_{v,e}\|_{L^2(e)}^2 \right),\label{eq:local estimator v}\\
	& \eta_{S}^2=h_e^2\|R_{\Gamma,S}\|_{L^2(S)}^2+	\left(\frac{1}{2} \sum_{\bm{x} \in \partial S}h_e \left|J_{\bm{x}}\right|^2\right)\label{eq:local estimator e}.
\end{align}
In particular, $\partial T \cap \Gamma$ contains at most two boundary edges.
Define 
\begin{equation*}
	\mathcal{R}_h(\Psi):=\ell(\Psi) -A(\bm{U}^h,\Psi), \quad \forall \Psi\in \mathcal{H}_{\Gamma}^1(\Omega;\Gamma).
\end{equation*}
It follows from \eqref{eq:variational problem} that
\begin{align*}
	\mathcal{R}_h(\Psi)&=-\rho_1\sum_{T\in \mathcal{M}_{h}}\left(\int_T\nabla w^h\cdot\nabla\overline{\psi}-k^2nw^h\overline{\psi} \right)\text{d}x
	+\rho_2\sum_{T\in \mathcal{M}_{h}}\left(\int_T\nabla v^h\cdot\nabla\overline{\phi}-k^2v^h\overline{\phi} \right)\text{d}x\\
	&\quad-\sum_{S\in \mathcal{S}_{h}}\left(\int_S \mu \nabla_\Gamma w^h\cdot\nabla_\Gamma\overline{\psi}+\lambda w^h\overline{\psi} \right)\text{d}s+\sum_{S\in \mathcal{S}_h}\int_S g\overline{\psi}\text{d}s.
\end{align*}
By using the integration by parts and \eqref{eq:integration by parts}, one can get
\begin{align*}
	\mathcal{R}_h(\Psi)&=\rho_1\sum_{T\in \mathcal{M}_{h}}\bigg(\int_T\left(\Delta w^h+k^2nw^h\right)\overline{\psi}\text{d}x-\sum_{e\in \partial T\setminus \Gamma}\int_e\nabla w^h \cdot \bm{n}\overline{\psi}\text{d}s\bigg)\\
	&\quad-\rho_2\sum_{T\in \mathcal{M}_{h}}\bigg(\int_T\left(\Delta v^h+k^2v^h\right)\overline{\phi}\text{d}x-\sum_{e\in \partial T\setminus \Gamma}\int_e\nabla v^h \cdot \bm{n}\overline{\phi}\text{d}s\bigg)\\
	&\quad+\sum_{S\in \mathcal{S}_{h}}\bigg(\int_S\left(g+\mu \Delta_\Gamma w^h-\lambda w^h-\rho_1\partial_{\bm \nu} w^h+\rho_2\partial_{\bm \nu} v^h\right)\overline{\psi}\text{d}s{-}\sum_{\bm{x}\in \partial S}\mu \nabla_{\Gamma}w^h(\bm{x}) \cdot \bm{n}_{\partial S}\overline{\psi}(\bm{x})\bigg).   
\end{align*}
which, by using the notations of residuals, can be written as
\begin{align}\label{eq:residual equation}
	\mathcal{R}_h(\Psi) &=\sum_{T\in \mathcal{M}_{h}}\bigg(\int_T R_{w,T}\overline{\psi}\text{d}x-\sum_{e\in \partial T\setminus \Gamma}\int_e\frac{1}{2}J_{w,e}\overline{\psi}\text{d}s\bigg)\nonumber\\
	&\quad-\sum_{T\in \mathcal{M}_{h}}\bigg(\int_TR_{v,T}\overline{\phi}\text{d}x-\sum_{e\in \partial T\setminus \Gamma}\int_e\frac{1}{2}J_{v,e}\overline{\phi}\text{d}s\bigg)\\
	&\quad+\sum_{S\in \mathcal{S}_{h}}\bigg(\int_SR_{\Gamma,S}\overline{\psi}\text{d}s{-}\sum_{\bm{x}\in \partial S}\frac{1}{2}J_{\bm{x}}\overline{\psi}(\bm{x})\bigg). \nonumber
\end{align}

To show the a posteriori error estimates, we use a uniform notation $\Pi_h$ to denote the Scott-Zhang interpolation in this subsection and denote $\Psi_h=(\psi_h,\phi_h)=\Pi_h \Psi$ for $\Psi=(\psi,\phi)\in \mathcal{H}_{\Gamma}^1(\Omega;\Gamma)$. It follows from the interpolation estimates that
\begin{align*}
	&\|\psi-\Pi_h \psi \|_{L^2(T)} \lesssim h_{T} \|\nabla \psi\|_{L^2(\widetilde{T})},\\
	&\|\psi-\Pi_h \psi \|_{L^2(e)} \lesssim h_{e}^{1/2} \| \psi\|_{H^1(\widetilde{T}_e)}.		
\end{align*}
Here, $\widetilde{T}$ and $\widetilde{T}_e$ are the union of all the elements in $\mathcal{M}_{h}$, which have nonempty intersection with the element $T$ and the side $e$, respectively. Besides, the above interpolation estimates can be extended to the boundary case such that
\begin{align*}
	\left\|\psi-\Pi_h\psi\right\|_{L^2\left(e\right)}  \lesssim h_e\left\|\nabla_{\Gamma} \psi\right\|_{L^2\left(\widetilde{e}\right)},\quad
	|\psi(\bm{x})-\Pi_h \psi(\bm{x}) | \lesssim h_{e}^{1/2} \| \nabla_{\Gamma} \psi\|_{L^2(\widetilde{\bm{x}})},
\end{align*}
where $\widetilde{e}$ and $\widetilde{\bm{x}}$ are the union of all the elements in $\mathcal{S}_{h}$ having nonempty intersection with the edge $e$ and point $\bm{x}$, respectively; see for example \cite{AFFKP15,AFKPP13,EM18}.

\begin{theorem}[Upper bound]
	Let $\bm{U}=(w,v)$ and $\bm{U}^h=(w^h,v^h)$ be the solutions of problems \eqref{eq:variational problem} and \eqref{eq:discrete problem}, respectively. Then the following a posteriori error estimate holds:
	\begin{equation*}
		\left\|\bm{U}-\bm{U}^h\right\|_{\mathcal{H}^1(\Omega; \Gamma)}
		\lesssim \left(\sum_{T\in \mathcal{M}_{h}}\eta_{T}^2\right)^{1/2} + \|f-f^h\|_{H^{1/2}(\Gamma)} .
	\end{equation*}    
\end{theorem}
\begin{proof}
	Set $\theta_f=f-f^h$, $\Theta_f=(0,\mathcal{L}\theta_f)$ by using the lifting function $\mathcal{L}$ defined in Lemma~\ref{lem:trace lifting}, then it is easy to check that
	\[\bm{E}_0:=(\bm U-\bm U^h)+\Theta_f \in \mathcal{H}_{\Gamma}^1(\Omega;\Gamma), \qquad \|\Theta_f\|_{\mathcal{H}^1(\Omega;\Gamma)} \lesssim \|f-f^h\|_{H^{1/2}(\Gamma)}, \]
	and
	\begin{equation}\label{ieq:upper bound inf-sup}
		\|\bm{E}_0\|_{\mathcal{H}^1(\Omega;\Gamma)}\lesssim \sup_{0\ne\Psi\in{\mathcal{H}^1_{\Gamma}(\Omega;\Gamma)}}
		\frac{|A(\bm{E}_0,\Psi)|}{\|\Psi\|_{\mathcal{H}^1(\Omega;\Gamma)}},   
	\end{equation}
	owing to the inf-sup condition \eqref{ieq:vpinf-sup}.
	For any $\Psi\in \mathcal{H}^1_\Gamma(\Omega;\Gamma)$, the following error representation holds
	\begin{align*}
		A(\bm E_0,\Psi)&=A(\bm{U},\Psi)-A(\bm{U}^h,\Psi)+A(\Theta_f,\Psi)\\
		&=A(\bm{U},\Psi-\Pi_h\Psi)-A(\bm{U}^h,\Psi-\Pi_h\Psi)+A(\Theta_f,\Psi)\\
		&=\mathcal{R}_h(\Psi-\Pi_h\Psi)+A(\Theta_f,\Psi).
	\end{align*}
	Applying the Scott-Zhang interpolation estimates to \eqref{eq:residual equation} gives
	\begin{equation}
		\mathcal{R}_h(\Psi-\Pi_h\Psi) \leq C\left(\sum_{T\in \mathcal{M}_{h}}\eta_{w,T}^2+\sum_{T\in \mathcal{M}_{h}}\eta_{v,T}^2+\sum_{S\in \mathcal{S}_{h}}\eta_{S}^2 \right)^{1/2}\|\Psi\|_{\mathcal{H}^1(\Omega; \Gamma)},
	\end{equation}
	and 
	\begin{equation}
		A(\Theta_f,\Psi)\leq C\|f-f^h\|_{H^{1/2}(\Gamma)}\|\Psi\|_{\mathcal{H}^1(\Omega; \Gamma)}.
	\end{equation}
	Noting that $\bm E=\bm E_0-\Theta_f$, the proof is completed by using \eqref{ieq:upper bound inf-sup} and the triangle inequality. 
\end{proof}

\subsection{A posteriori error estimates: local lower bound}
\label{sec:3.4}

In this subsection, we further establish a local lower bound to ensure the efficiency of the local error estimators. 

For any $V=T,\partial T, S$ with $T\in \mathcal{M}_{h}$ and $S\in \mathcal{S}_{h}$, we define the orthogonal projection $Q_V^p:L^2(V)\rightarrow P_p(V)$ through
\[(Q_V^p\psi,\phi)_V=(\psi,\phi)_V,\quad \forall \phi \in P_p(V),\]where $(\cdot,\cdot)_V$ denotes the inner product on $L^2(V)$. Denote
\begin{align*}
	&\widehat{R}_{w,T}=Q_T^1R_{w,T},\qquad 
	\widehat{R}_{v,T}=Q_T^1R_{v,T},\qquad \widehat{R}_{\Gamma,S}=Q_S^1R_{\Gamma,S}, \\ &\widehat{J}_{w,e}=Q_{\partial T}^0 J_{w,e},\qquad 
	\widehat{J}_{v,e}=Q_{\partial T}^0 J_{v,e}.
\end{align*}
Obviously, $\widehat{R}_{v,T}=R_{v,T}, \widehat{J}_{w,e}=J_{w,e},  
\widehat{J}_{v,e}=J_{v,e}$.

\begin{theorem}[Local lower bound]
	Let $\bm{U}=(w,v)$ and $\bm{U}^h=(w^h,v^h)$ be the solutions of \eqref{eq:variational problem} and \eqref{eq:discrete problem}, respectively. There exists a constant $h_0>0$ such that for all $0<h<h_0$ and any $T\in \mathcal{M}_h$ the following estimate holds:
	\begin{equation*}
		\eta_{T}^2
		\lesssim \|\bm U-\bm U^h\|^2_{\mathcal{H}^1(T^*;\Gamma_{T^*})} + osc^2_{T^*},
	\end{equation*} 
	where $T^*$ is the union of all elements sharing at least one common side with $T$ or one common vertex located on $\Gamma$ with $T$, $\Gamma_{T^*}=\Gamma \cap \overline{T^*}$, and
	\begin{equation*}
		osc^2_{T^*}=h_{T^{'}}^2\sum_{T^{'}\subset T^*} \|  R_{w,T^{'}}- \widehat{R}_{w,T^{'}} \|_{L^2(T^{'})}^2+h_e^2\sum_{S^{'}\subset{\Gamma_{T^*}} }\| R_{\Gamma,S^{'}}-\widehat{R}_{\Gamma,S^{'}}\| _{L^2(S^{'})}^2.
	\end{equation*}
\end{theorem}

\begin{proof}
	From the definition of $\mathcal{R}_h(\Psi)$, for any $\Psi \in \mathcal{H}_{\Gamma}^1(\Omega,\Gamma)$,
	\begin{align}\label{eq:residual equation 2}
		A(\bm U-\bm U^h,\Psi) = \mathcal{R}_h(\Psi) &=\sum_{T\in \mathcal{M}_{h}}\bigg(\int_T R_{w,T}\overline{\psi}\text{d}x-\int_TR_{v,T}\overline{\phi}\text{d}x\bigg)\nonumber\\
		&\quad-\sum_{e\in \partial T\setminus \Gamma}\bigg(\int_eJ_{w,e}\overline{\psi}\text{d}s-\int_eJ_{v,e}\overline{\phi}\text{d}s\bigg)\\
		&\quad+\sum_{S\in \mathcal{S}_{h}}\int_SR_{\Gamma,S}\overline{\psi}\text{d}s-\sum_{\bm{x}\in \partial S}J_{\bm{x}}\overline{\psi}(\bm{x}). \nonumber
	\end{align}
	(i) For any element $T\in \mathcal{M}_h$ ,  let $\psi_T=b_T\widehat{R}_{w,T} \in H_0^1(T)$ where $b_T=27\prod_{j=1}^{3} \lambda_j^T$ is a canonical bubble function satisfying $0\leq b_T\leq1,b_T|_{\partial T}=0,\mathrm{supp}(\psi_T)\subset \overline{T}$ with $\lambda_j^T$ being the barycentric coordinate functions on $T$. By taking $\psi=\psi_T$ and $\phi=0$ in \eqref{eq:residual equation 2}, one obtains
	\begin{align}\label{eq:Rwt equaility}
		\int_TR_{w,T}\overline{\psi_T}\text{d} x&=\rho_1\int_T(\nabla (e_w)\cdot \nabla \overline{\psi_T}-k^2n(e_w)\overline{\psi_T})\text{d} x.
	\end{align}
	Furthermore,
	\begin{equation*}
		\|\psi_T\|_{L^2(T)} \lesssim\|\widehat R_{w,T}\|_{L^2(T)},\quad 
		\|\nabla\psi_T\|_{L^2(T)} \lesssim h_T^{-1} \|\widehat R_{w,T}\|_{L^2(T)}.
	\end{equation*}
	Therefore, by using the Cauchy–Schwarz inequality and the estimates for the bubble function then yields 
	\begin{align*}
		\int_TR_{w,T}\overline{\psi_T}\text{d} x \lesssim\|w-w^h\|_{H^1(T)}\|\psi_T\|_{H^1(T)}\lesssim h_T^{-1}\|w-w^h\|_{H^1(T)}\|\widehat{R}_{w,T}\|_{L^2(T)}.
	\end{align*}
	Since ${\widehat{R}_{w,T}}$ is a polynomial function in $T$, we have 
	\begin{align*}
		\|\widehat{R}_{w,T}\|^2_{L^2(T)}\lesssim\int_T|\widehat{R}_{w,T}|^2\overline{b_T}\text{d} x &=\int_T\widehat{R}_{w,T}\overline{\psi_T} \text{d} x\\
		&\leq\int_T(\widehat{R}_{w,T}-R_{w,T})\overline{\psi_T} \text{d} x+\int_TR_{w,T}\overline{\psi_T} \text{d} x\\
		&\lesssim\|\widehat{R}_{w,T}-R_{w,T}\|_{L^2(T)}\|\widehat{R}_{w,T}\|_{L^2(T)}+\int_TR_{w,T}\overline{\psi_T} \text{d} x,
	\end{align*}
	which further implies
	\begin{equation*}
		\|\widehat{R}_{w,T}\|^2_{L^2(T)}\lesssim \left(\|\widehat{R}_{w,T}-R_{w,T}\|_{L^2(T)}+ C h_T^{-1}\|w-w^h\|_{H^1(T)}\right)\|\widehat{R}_{w,T}\|_{L^2(T)}.
	\end{equation*}
	This, together with the triangle inequality, gives
	\begin{align}\label{ieq:apea lower bound Rwt}
		h^2_T\|R_{w,T}\|^2_{L^2(T)}
		&\lesssim h^2_T\|\widehat{R}_{w,T}\|^2_{L^2(T)}+h^2_T\|R_{w,T}-\widehat{R}_{w,T}\|^2_{L^2(T)} \nonumber\\
		&\lesssim\|w-w^h\|^2_{H^1(T)}+h^2_T\|\widehat{R}_{w,T}-R_{w,T}\|^2_{L^2(T)}.
	\end{align}
	On the other hand, since $\widehat{R}_{v,T}=R_{v,T}$, taking the test function $\psi = 0, \phi = b_T{R}_{v,T}$ in \eqref{eq:residual equation 2}, the same argument leads to
	\begin{equation}\label{ieq:apea lower bound Rvt}
		h^2_T\|R_{v,T}\|^2_{L^2(T)} \lesssim\|v-v^h\|^2_{H^1(T)}.
	\end{equation}
	(ii) For any interior edge $e=T_1\cap T_2$, let $\psi_e=b_e {J}_{w,e}\in H_0^1(\Omega_e)$ where $\Omega_e=T_1\cup T_2, b_e=4\prod_{j=1}^{2} \lambda_j^e$ is an edge bubble function and $\lambda_j^e$ is the barycentric coordinate functions associated with the nodes of $e$. The properties of the edge bubble function imply that
	\begin{equation*}
		\|\psi_e\|_{L^2(\Omega_e)} \lesssim h_e^{1/2}\|J_{w,e}\|_{L^2(e)},\quad
		\|\psi_e\||_{H^1(\Omega_e)} \lesssim h_e^{-1/2} \|J_{w,e}\|_{L^2(e)}.
	\end{equation*}
	Thus, taking  the test function $(\psi_e, 0)$ in \eqref{eq:residual equation 2} and using integration by parts give
	\begin{align*}
		\int_eJ_{w,e}\overline{\psi_e}\text{d} x &=\sum_{T\subset{\Omega_e}}\int_TR_{w,T}\overline{\psi_e}\text{d}x -\rho_1\int_{\Omega_e}\left(\nabla (w-w^h)\cdot\nabla \overline{\psi_e}-k^2n(w-w^h)\overline{\psi_e}\right)\text{d} x\\
		&\lesssim \|w-w^h\|_{H^1(\Omega_e)}\|\psi_e\|_{H^1(\Omega_e)}+\sum_{T\subset\Omega_e}\left(\|R_{w,T}\|_{L^2(T)}\|\psi_e\|_{L^2(T)}\right)\\
		&\lesssim h_e^{-1/2}\|w-w^h\|_{H^1(\Omega_e)}\|{J}_{w,e}\|_{L^2(e)}+ h_e^{1/2}\sum_{T\subset\Omega_e}\|R_{w,T}\|_{L^2(T)}\|{J}_{w,e}\|_{L^2(e)}.
	\end{align*}
	In addition,
	\begin{align*}        
		\|{J}_{w,e}\|^2_{L^2(e)} \lesssim\int_e|{J}_{w,e}|^2\overline{b_e}\text{d} s         &=\int_e{J}_{w,e}\overline{\psi_e} \text{d} s.
	\end{align*}
	Therefore,
	\begin{equation*}
		\|J_{w,e}\|_{L^2(e)}\lesssim h_e^{-1/2} \|w-w^h|_{H^1(\Omega_e)}+h_e^{1/2}\sum_{T\subset\Omega_e}\|R_{w,T}\|_{L^2(T)}.
	\end{equation*}
	Multiplying by $h_e^{1/2}$ and using the shape regularity of the mesh, which implies $h_e\simeq h_T$ for $T\subset\Omega_e$, yield
	\begin{equation}
		h_e\|J_{w,e}\|_{L^2(e)}^2 \lesssim \|w-w^h\|_{H^1(\omega_e)}^2 + \sum_{T\subset\omega_e}
		h_T^2 \|R_{w,T}\|_{L^2(T)}^2 .
	\end{equation}
	Analogously,
	\begin{align}\label{ieq:apea lower bound Jve}
		h_e\|J_{v,e}\|^2_{L^2(e)}
		\lesssim \|v-v^h\|^2_{H^1(\Omega_e)}+\sum_{T\subset\Omega_e}h_T^2\|R_{v,T}\|_{L^2(T)}^2.
	\end{align}
	(iii) For any $S\in \mathcal{S}_h$ which is one face of $\Omega_S\in\mathcal{M}_h$, set $\psi _S=b_S\widehat{R}_{\Gamma,S}$ where $b_S$ is also a standard edge bubble function defined as $b_e$. Analogous to case (ii), we have 
	\begin{equation*}
		\|\psi_S\|_{L^2(\Omega_S)} \lesssim h_e^{1/2}\|\widehat R_{\Gamma,S}\|_{L^2(S)},\quad 
		\|\psi_S\||_{H^1(\Omega_S)} \lesssim h_e^{-1/2} \|\widehat R_{\Gamma,S}\|_{L^2(S)}.
	\end{equation*}
	Since $\psi _S=b_S\widehat{R}_{\Gamma,S}$ is a polynomial of degree 2 on $S$, the standard inverse estimate gives
	\begin{equation*}
		\|\psi_S\|_{L^2(S)} \lesssim \|\widehat R_{\Gamma,S}\|_{L^2(S)},\quad 
		\|\nabla_\Gamma\psi_S\|_{L^2(S)} \lesssim h_e^{-1}\|\widehat R_{\Gamma,S}\|_{L^2(S)}.
	\end{equation*}
	Hence, 
	\begin{equation*}
		\|\psi_S\|_{H^1(S)} \lesssim h_e^{-1}\|\widehat R_{\Gamma,S}\|_{L^2(S)}.
	\end{equation*}  
	By taking $\Psi=(\psi_S,\psi_S)$ in \eqref{eq:residual equation 2}, we get
	\begin{align*}
		\int_SR_{\Gamma,S}\overline{\psi_S}\text{d} s &=\rho_1\int_{\Omega_S}\left(\nabla (w-w^h)\cdot\nabla \overline{\psi_S}-k^2n(w-w^h)\overline{\psi_S}\right)\text{d} x-\rho_2\int_{\Omega_S}\left(\nabla (v-v^h)\cdot\nabla \overline{\psi_S}-k^2(v-v^h)\overline{\psi_S}\right)\text{d} x
		\\& \quad +\int_S\left(\mu \nabla_\Gamma (w-w^h)\cdot\nabla_\Gamma\overline{\psi_S}+\lambda(w-w^h)\overline{\psi_S}\right)\text{d} s -\sum_{T\subset{\Omega_S}}\left(\int_{T}R_{w,T}\overline{\psi_S}\text{d}x- \int_{T}R_{v,T}\overline{\psi_S}\text{d}x \right)\\
		&\lesssim h_e^{-1/2}\left(\|w-w^h\|_{H^1(\Omega_S)}+\|v-v^h\|_{H^1(\Omega_S)}\right)\|\widehat{R}_{\Gamma,S}\|_{L^2(S)}+h_e^{-1}\|w-w^h\|_{H^1(S)}\|\widehat{R}_{\Gamma,S}\|_{L^2(S)}\\    &\quad+h_e^{1/2}\sum_{T\subset\Omega_S}\left(\|R_{w,T}\|_{L^2(T)}+\|R_{v,T}\|_{L^2(T)}\right)\|\widehat{R}_{\Gamma,S}\|_{L^2(S)}.
	\end{align*}
	Hence,
	\begin{align}\label{ieq:apea lower bound Rs}
		h^2_e\|R_{\Gamma,S}\|^2_{L^2(S)} &\lesssim\|w-w^h\|^2_{H^1(\Omega_S)}+\|v-v^h\|^2_{H^1(\Omega_S)}+\|w-w^h\|^2_{H^1(S)} \nonumber\\   &\quad+\sum_{T\subset\Omega_S}h^2_T\left(\|R_{w,T}\|^2_{L^2(T)}+\|R_{v,T}\|^2_{L^2(T)}\right)+h_e^2\|\widehat{R}_{\Gamma,S}-R_{\Gamma,S}\|^2_{L^2(S)}.
	\end{align}
	
	(iv) For any endpoint ${\bm{x}} =\partial S_1\cap \partial S_2, S_1,S_2\in\mathcal{S}_h$, define $\Omega_{\bm x}=\bigcup\{ T\in \mathcal{M}_h: \bm{x}\in \overline{T} \}$ and $S_{\bm{x}}=S_1\cup S_2$. Let $\widehat{\psi} _{\bm x}=b_{\bm{x}}J_{\bm{x}}$, where $b_{\bm{x}}$ is a hat function with $b_{\bm{x}}(\bm x)=1$. Additionally, let $\psi_{\bm{x}}=\varphi_{\bm{x}}J_{\bm{x}}$ where $\varphi_{\bm{x}}$ is the standard $P_1$ nodal basis function, which is supported on $\overline{\Omega_{\bm{x}}}$, associated with the vertex $\bm{x}$. In particular, $\varphi_{\bm{x}}|_{S_{\bm{x}}}=b_{\bm{x}}$. Therefore, $\psi_{\bm{x}}|_{S_{\bm{x}}}=\widehat{\psi} _{\bm x}$
	and $\psi_{\bm{x}}$ vanishes on $\Gamma \setminus S_{\bm{x}}$. Note that the number of elements in $\Omega_{\bm{x}}$ is uniformly bounded. Hence,
	\begin{align*}
		\|\psi_{\bm{x}}\|_{L^2({\Omega_{\bm{x}}})}\lesssim h_e|J_{\bm{x}}|,\quad
		\|\psi_{\bm{x}}\|_{H^1({\Omega_{\bm{x}}})}\lesssim |J_{\bm{x}}|,
	\end{align*}
	and
	\begin{align*}
		\|\widehat{\psi} _{\bm x}\|_{L^2({S_{\bm{x}}})}\lesssim h_e^{1/2}|J_{\bm{x}}|,\quad
		\|\widehat{\psi} _{\bm x}\|_{H^1({S_{\bm{x}}})}\lesssim h_e^{-1/2}|J_{\bm{x}}|.
	\end{align*}
	For those edges $e\in \bigcup_{T\in\Omega_{\bm x}}\partial T\backslash \mathcal{S}_h$, we have
	\begin{equation*}
		\|\psi_{\bm{x}}\|_{L^2({e})}\lesssim h_e^{1/2}|J_{\bm{x}}|.
	\end{equation*}
	Then taking $\Psi=\left(\psi_{\bm{x}},\psi_{\bm{x}} \right)$ in \eqref{eq:residual equation 2} yields
	\begin{align*}
		|J_{\bm{x}}|^2= J_{\bm{x}} \overline{\psi_{\bm{x}}(\bm{x})} 
		&=-A(\bm E,\Psi_{\bm x})+\sum_{T\subset{\Omega_{\bm{x}}}} \left(\int_TR_{w,T}\overline{\psi_{\bm{x}}}\text{d}x-\int_TR_{v,T}\overline{\psi_{\bm{x}}}\text{d}x\right)\\
		&\quad-\sum_{e\in\cup_{T\in\Omega_{\bm x}}\partial T\backslash \mathcal{S}_h}\left(\int_eJ_{w,e}\overline{\psi_{\bm{x}}}\text{d}s-\int_eJ_{v,e}\overline{\psi_{\bm{x}}}\text{d}s\right) +\sum_{S'\subset{S_{\bm{x}}}}\int_{S'} R_{\Gamma,S'}\overline{\widehat{\psi}_{\bm{x}}}\text{d}s . 
	\end{align*}
	It is easy to estimate that
	\begin{equation*}    
		A(\bm E,\Psi_{\bm x})\lesssim \left(\|w-w^h\|_{H^1(\Omega_{\bm{x}})}+\|v-v^h\|_{H^1(\Omega_{\bm{x}})}+h_e^{-1/2}\|w-w^h\|_{H^1(S_{\bm{x}})}\right)|J_{\bm{x}}|.
	\end{equation*}
	For any $S'\subset S_{\bm x}$,
	\begin{equation*}     \int_{S'}R_{\Gamma,S'}\overline{\widehat{\psi}_{\bm x}}\text{d} s\lesssim h_e^{1/2}\|R_{\Gamma,S'}\|_{L^2(S')}|J_{\bm{x}}|.
	\end{equation*}
	Therefore,
	\begin{align}\label{ieq:apea lower bound Jx}
		h_e|J_{\bm{x}}|^2 &\lesssim \|w-w^h\|^2_{H^1(\Omega_{\bm{x}})}+\|v-v^h\|^2_{H^1(\Omega_{\bm{x}})}+\|w-w^h\|^2_{H^1(S_{\bm{x}})} \nonumber\\   &\quad+\sum_{T\subset\Omega_{\bm{x}}}h^2_T\left(\|R_{w,T}\|_{L^2(T)}^2+\|R_{v,T}\|^2_{L^2(T)}\right)+\sum_{e\in\cup_{T\in\Omega_{\bm x}}\partial T\backslash \mathcal{S}_h}h_e\left(\|J_{w,e}\|^2_{L^2(e)}+\|J_{v,e}\|^2_{L^2(e)}\right)\nonumber\\    &\quad+\sum_{S'\subset{S_{\bm{x}}}}h^2_e\|R_{\Gamma,S'}\|^2_{L^2(S')}.
	\end{align}
	Combining \eqref{ieq:apea lower bound Rwt}-\eqref{ieq:apea lower bound Jx} completes the proof.
\end{proof}

\section{Numerical experiments}
\label{sec:4}

In this section, several numerical examples will be presented to verify the accuracy and efficiency of the proposed adaptive FEM for solving the interior transmission problem with GTBCs as well as its extension to the calculation of interior transmission eigenvalues. The unit disk, an L-shaped domain $\Omega=(-1,1)\times(-1,1)\backslash[0,1)\times(-1,0]$  and a U-shaped domain $\Omega = (-1, 1)\times(-1, 1)\backslash (-1, 0.6]\times[-0.5, 0.5]$ will be considered to illustrate the performance of the method. All numerical results presented in this paper were obtained using MATLAB implementations of the algorithms on a laptop with CPU Core(TM) i5-1240P.

Define an effectivity index as follows \begin{equation*}
	I_{\mathrm{eff}}=\frac{\epsilon_h}{\left\|\bm{U}-\bm{U}^h\right\|_{\mathcal{H}^1(\Omega; \Gamma)}}\quad \mbox{with}\quad \epsilon_h = \left(\sum_{T\in \mathcal{M}_h} \eta_{T}^2\right)^{1/2}.
\end{equation*}

\begin{figure}
	\centering
	\begin{tabular}{ccc}
		\includegraphics[scale=0.27]{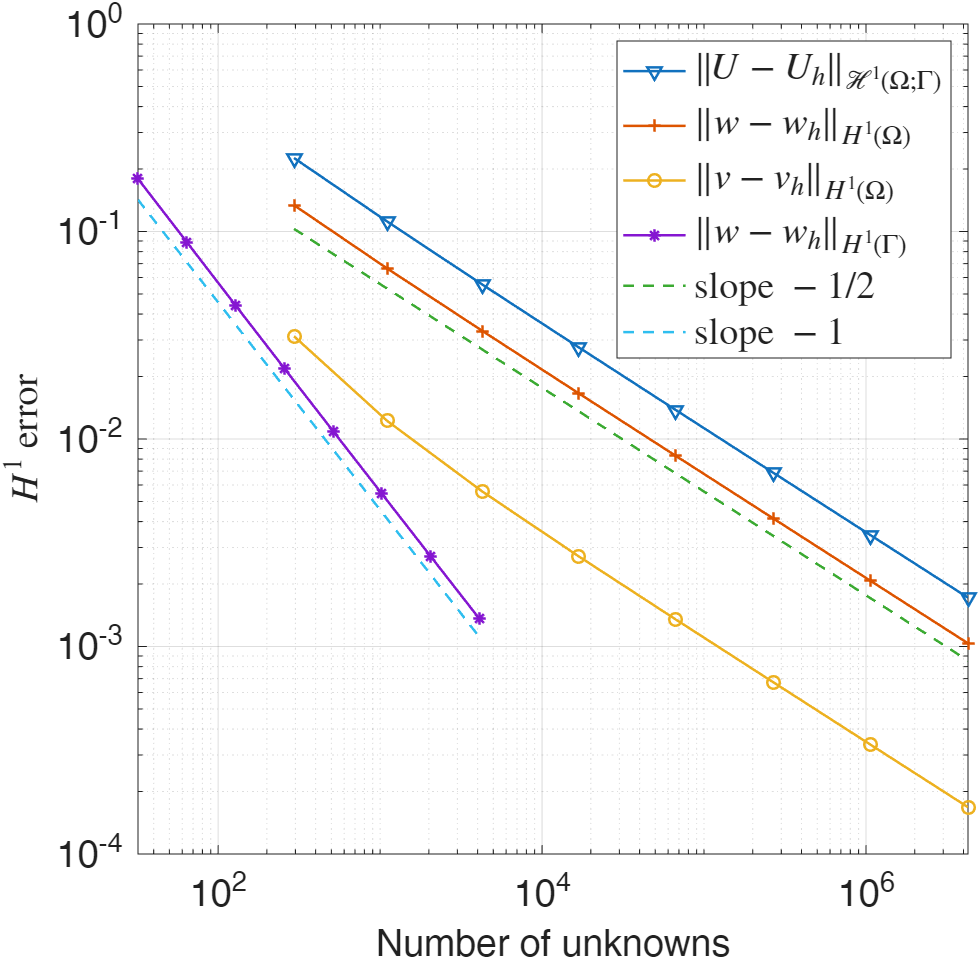} & 
		\includegraphics[scale=0.27]{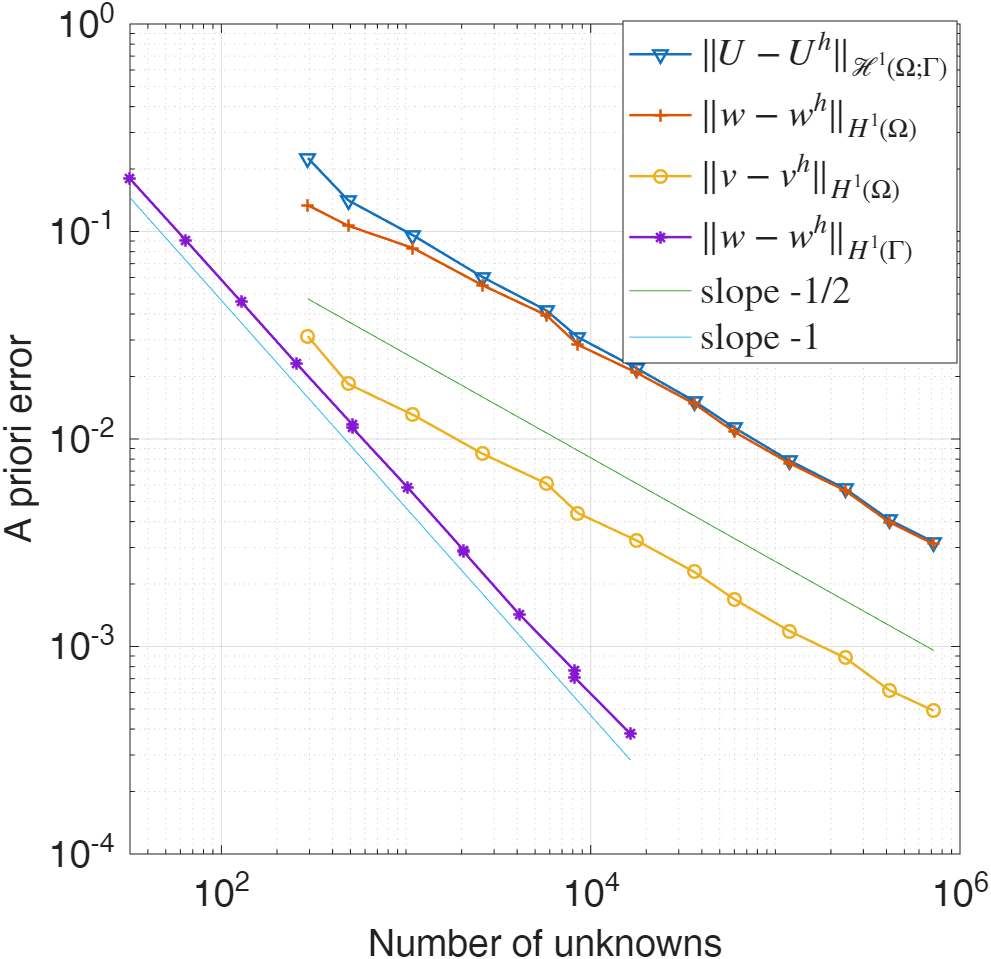} & 
		\includegraphics[scale=0.27]{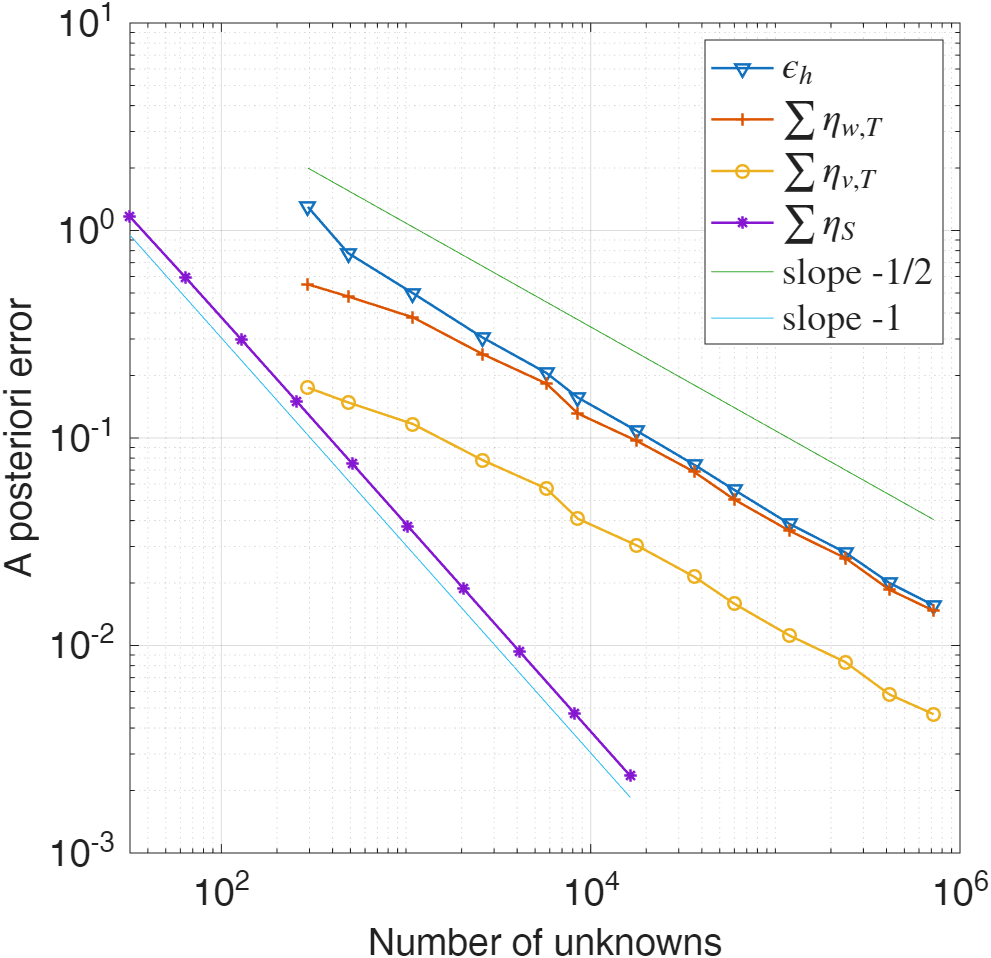}\\
		(a)  & (b)  & (c) 
	\end{tabular}
	\caption{Example 1. Convergence of the $H^1$-errors under uniform mesh refinement (a) and the $H^1$-errors, the a posteriori error estimators under adaptive refinement (b,c) for the problem in a unit disk.}
	\label{ex1_error}
\end{figure}

\textbf{Example 1} ({\it unit disk}). In this example, let $\Omega$ be the unit disk and set $k=1, n=4, m=1, \rho_1=1, \rho_2=2, \mu=3, \lambda=0.5$. Then the functions $w(r,\theta)=J_m(2r)e^{im\theta}$ and $v(r,\theta)=J_m(r)e^{im\theta}$ in polar coordinates can be selected as the exact solutions of the interior transmission problem \eqref{eq:Strong Model} with the boundary data $f=\left(J_m(2r)-J_m(r)\right)e^{im\theta}$ and $g=\rho_1\partial_rw-\rho_2\partial_rv+(\mu m^2+\lambda)J_m(2r)e^{im\theta}$. Here, $J_m$ is the first kind Bessel function of order $m$ and 
\begin{align*}
	\partial_rw&= 2\left(J_{m-1}(2r)-\frac{m}{2r}J_m(2r)\right)e^{im\theta} ,\\
	\partial_rv&=\left(J_{m-1}(r)-\frac{m}{r}J_m(r)\right)e^{im\theta}.
\end{align*}
Under the uniform refinement of finite element mesh, the convergence of $H^1$-errors with respect to the number of unknowns $\text{DoF}_h$ shown in Figure~\ref{ex1_error}(a) clearly verify the optimal convergence order $\mathcal{O}\left(\mathrm{DoF}_h^{-1 / 2}\right)$ in $\Omega$ and $\mathcal{O}\left(\mathrm{DoF}_h^{-1}\right)$ on $\Gamma$ which can also be observed under the adaptive refinement, see Figure~\ref{ex1_error}(b,c). 


\begin{figure}
	\centering
	\begin{tabular}{ccc}
		\includegraphics[scale=0.27]{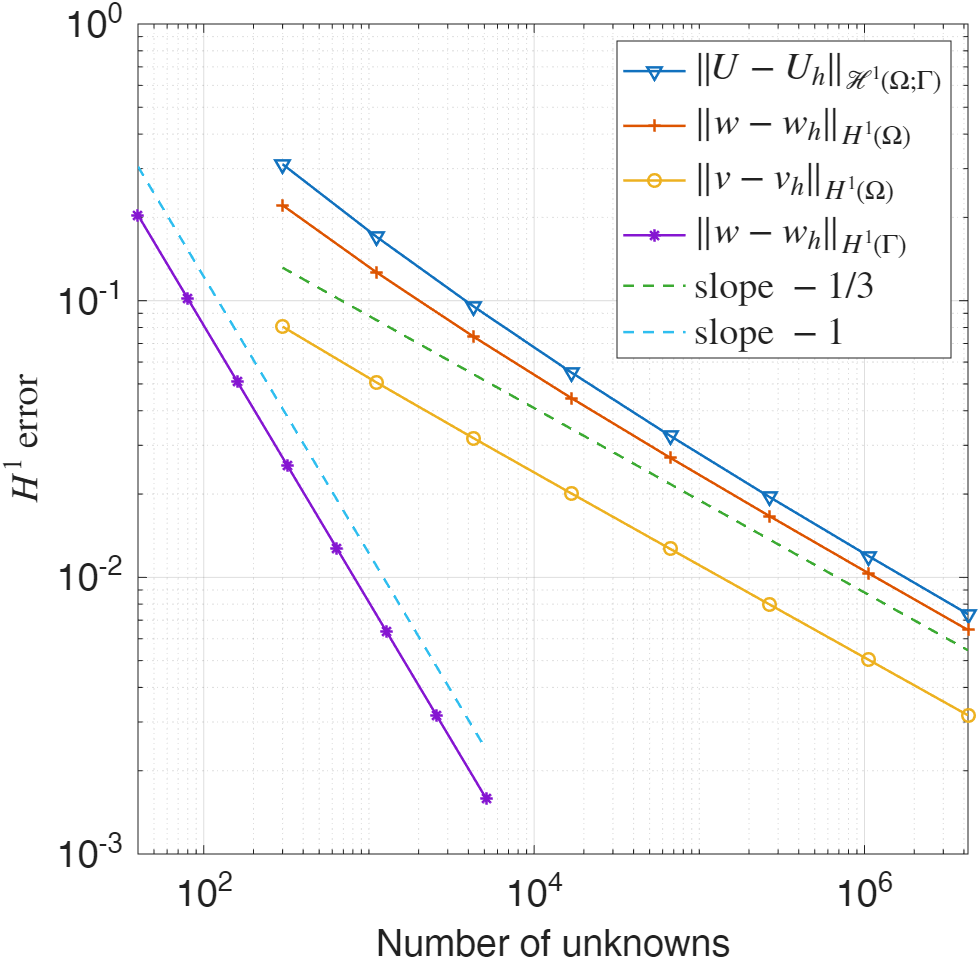} & 
		\includegraphics[scale=0.27]{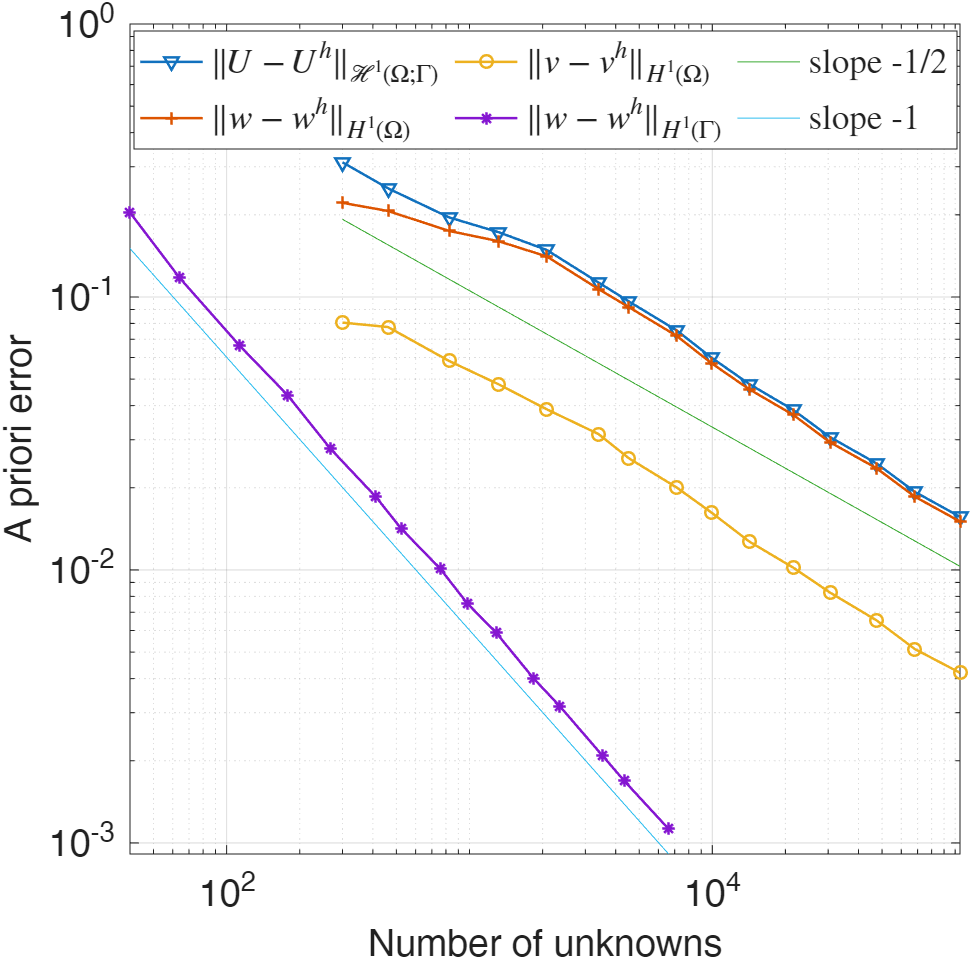} & 
		\includegraphics[scale=0.27]{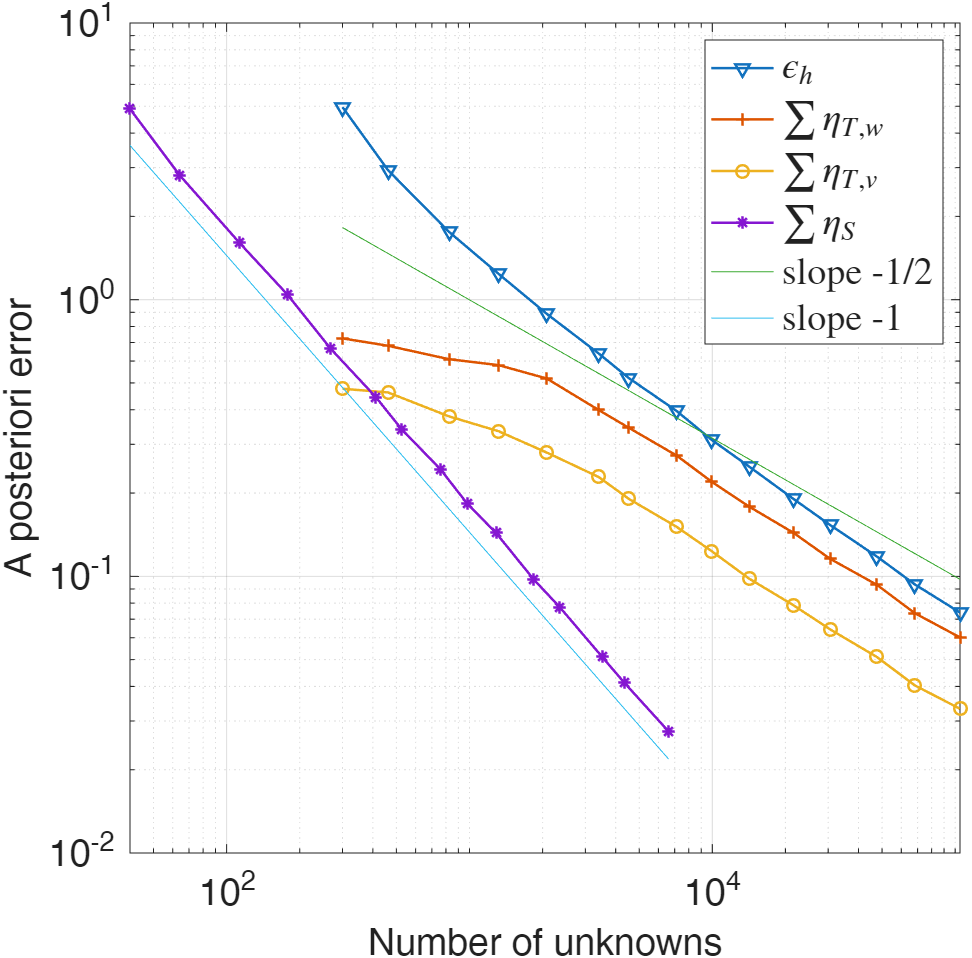}\\
		(a)  & (b)  & (c) 
	\end{tabular}
	\caption{Example 2. Convergence of the $H^1$-errors under uniform mesh refinement (a) and the $H^1$-errors, the a posteriori error estimators under adaptive refinement (b,c) for the problem in an L-shaped domain.}
	\label{ex2:error}
\end{figure}

\begin{figure}
	\centering   
	\begin{tabular}{ccc}
		\includegraphics[scale=0.27]{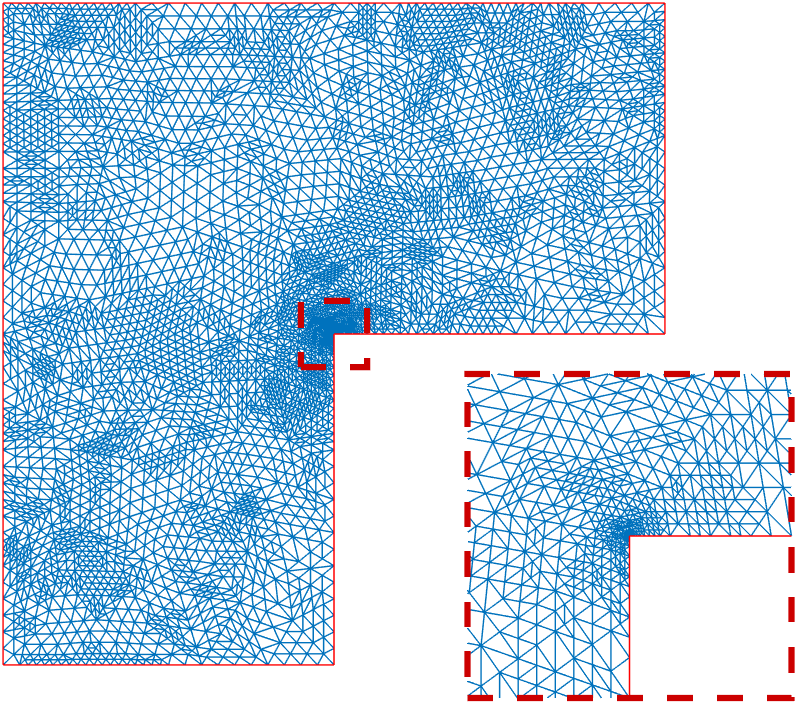} & 
		\includegraphics[scale=0.27]{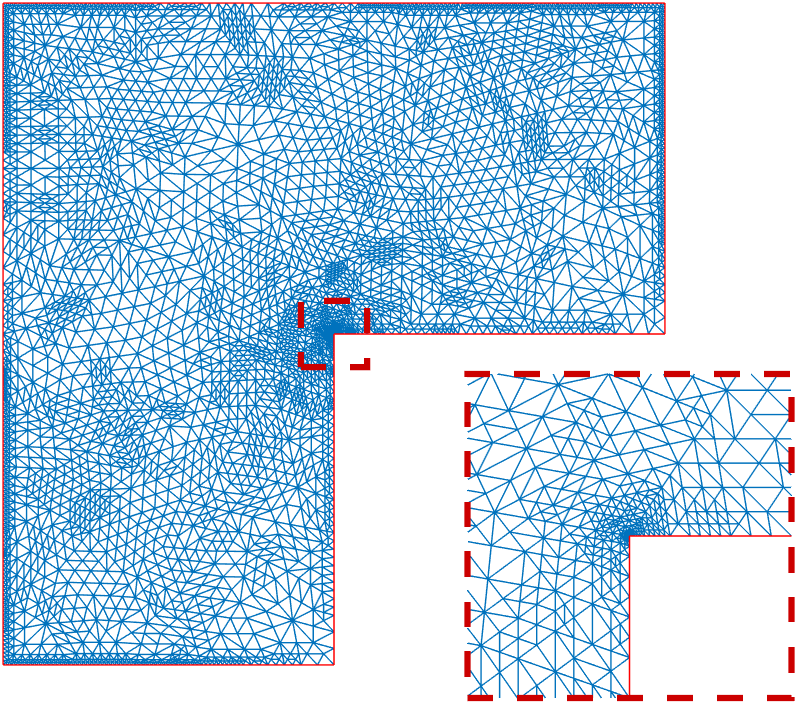} & 
		\includegraphics[scale=0.27]{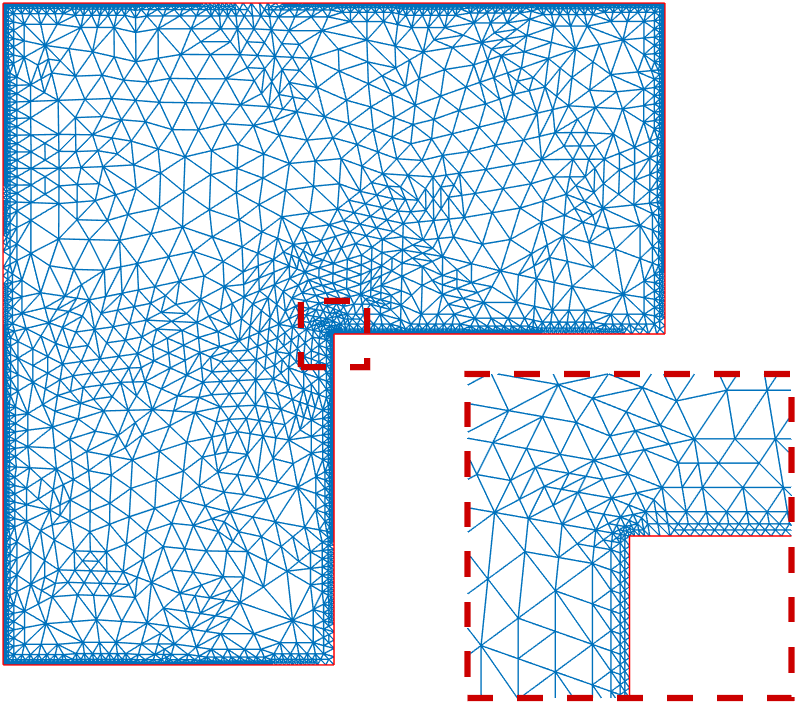}\\
		(a) $\mu=0.1, \text{Dof}_h=9470$  & (b) $\mu=1, \text{Dof}_h=8518$ & (c) $\mu=10, \text{Dof}_h=9658$
	\end{tabular}
	\caption{Example 2. Adaptive meshes for different values of $\mu$ for the problem in an L-shaped domain.}
	\label{ex2:meshL comparision}
\end{figure}

\begin{figure}
	\centering    
	\begin{tabular}{ccc}
		\includegraphics[scale=0.27]{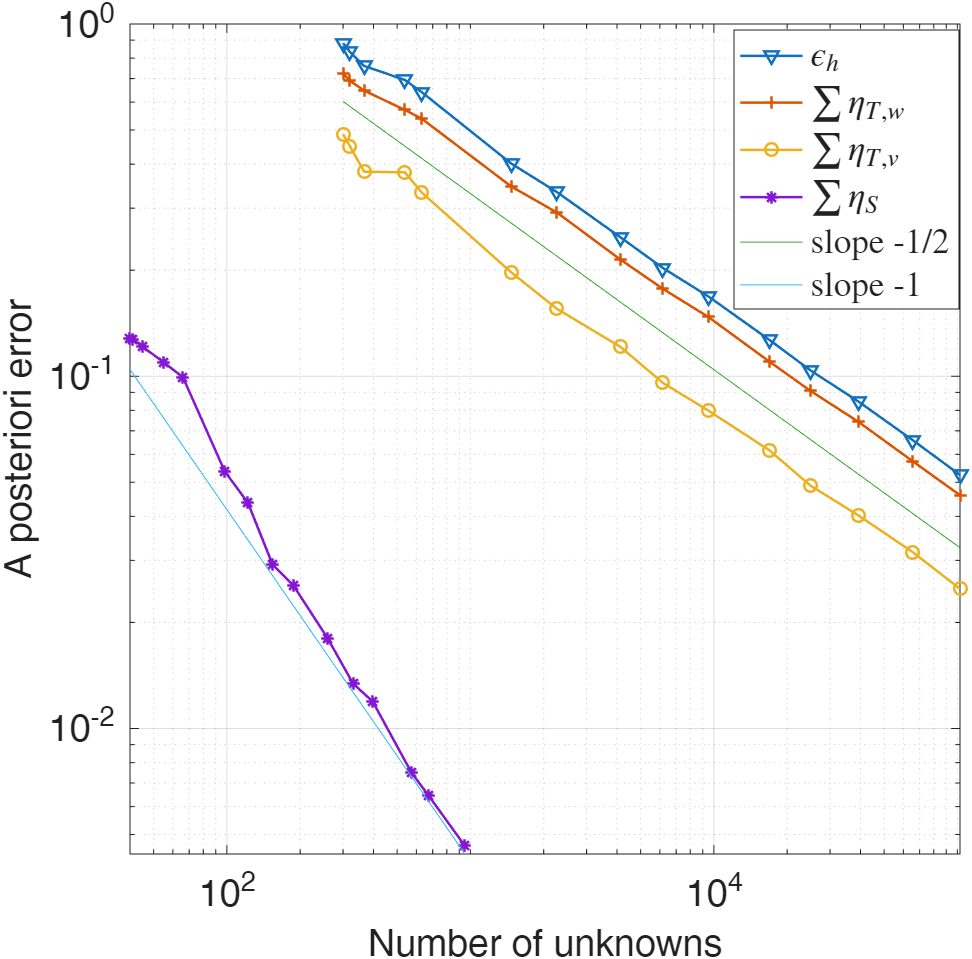} & 
		\includegraphics[scale=0.27]{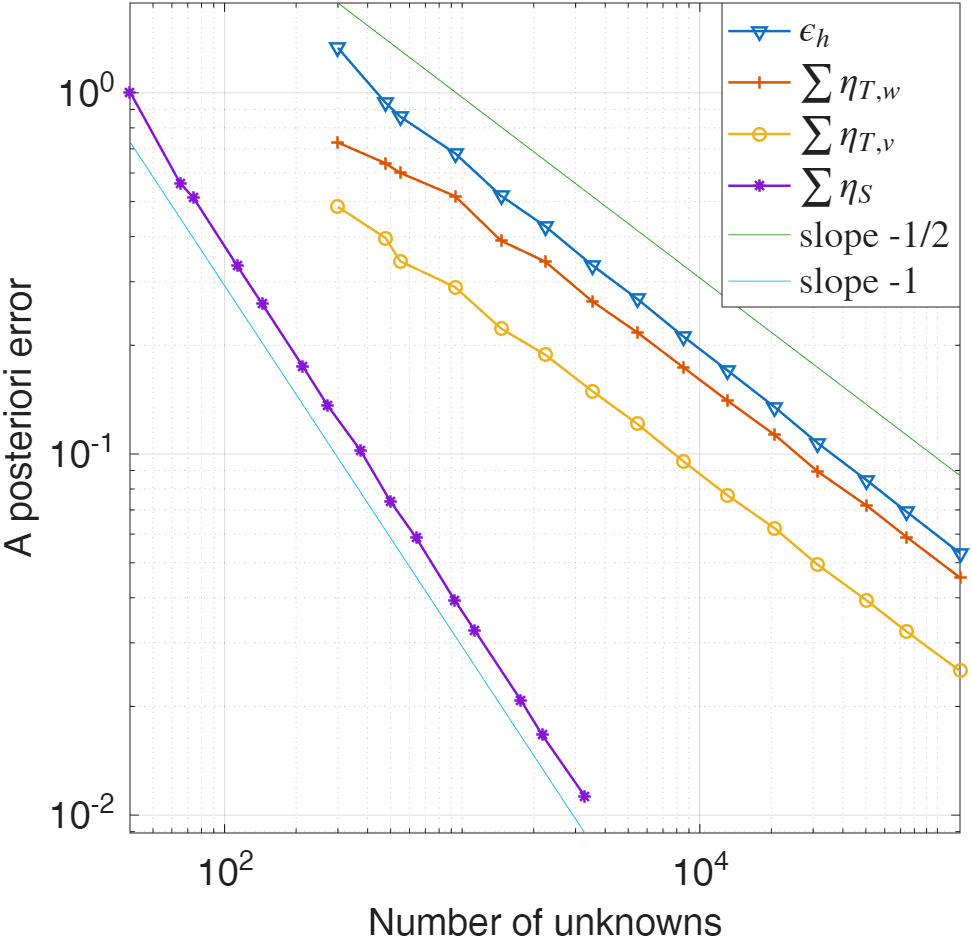} & 
		\includegraphics[scale=0.27]{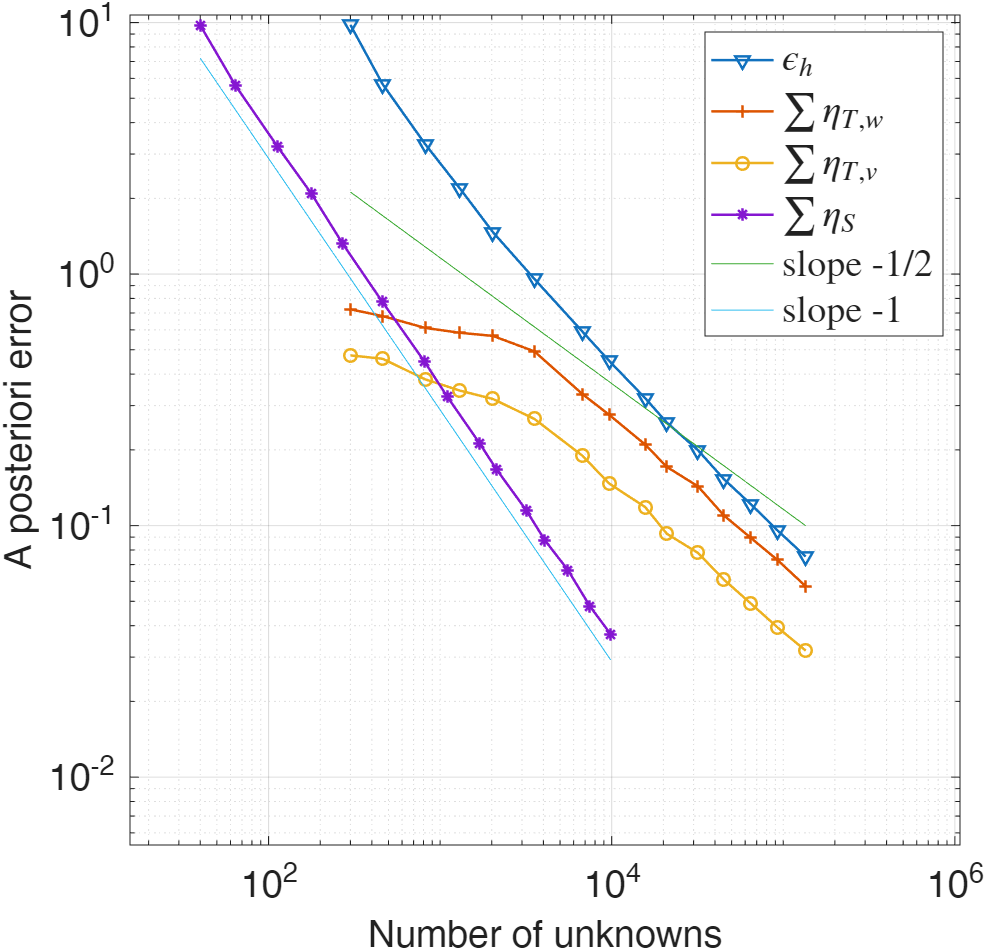}\\
		(a) $\mu=0.1$  & (b) $\mu=1$ & (c) $\mu=10$
	\end{tabular}
	\caption{Example 2. Convergence of the a posteriori error estimators for different values of $\mu$.}
	\label{ex2:error comparision}
\end{figure}

\textbf{Example 2} ({\it L-shaped domain}). In this example, let $\Omega$ be the L-shaped domain and set $k=1, n=4, \rho_1=1, \rho_2=2, \mu=5, \lambda=1.5-2\mathrm{i}$. We choose appropriate $f,g$ such that the exact solutions of the interior transmission problem \eqref{eq:Strong Model}, which admit a certain corner singularity near the point $(0,0)$
in this case, are given by
\begin{equation*}
	w(x_1,x_2)=J_{2/3}(2r)\sin\frac{2\theta}{3}+F(x_1,x_2),\quad 
	v(x_1,x_2)=2^{-1/3}J_{2/3}(r)\sin{\frac{2\theta}{3}},
\end{equation*}
where $F(x_1,x_2)$ is a correction function to ensure the compatibility condition \eqref{{eq:d}} and takes the form 
\begin{align*}
	F(x_1,x_2) = & A[\cos(2x_1)+\cos(2x_2)]+B[\sin(2x_1)-\sin(2x_2)]\\
	&+C\cos\left(\sqrt{2}(x_1-x_2)\right)+D\sin\left(\sqrt{2}(x_1-x_2)\right),   
\end{align*}
with 
\begin{align*}
	A&=-0.1744072606513751,\qquad B=-0.2672429488956379,\\
	C&=-0.0714070902529113,\qquad D=0.3779386027767911.
\end{align*}
The asymptotic property of the Bessel function implies that $w,v\in H^{5/3-\epsilon}(\Omega)$ for any $\epsilon>0$. As a result, a uniform mesh refinement only gives the reduced convergence order $\mathcal{O}\left(\mathrm{DoF}_h^{-1/3}\right)$ in $\Omega$ which is in agreement with the predicted regularity-limited rates, but optimal order $\mathcal{O}\left(\mathrm{DoF}_h^{-1}\right)$ on $\Gamma$, see Figure~\ref{ex2:error}(a). In contrast, the adaptivity substantially improves the convergence behavior and ensures the quasi-optimal convergence $\mathcal{O}\left(\mathrm{DoF}_h^{-1/2}\right)$ in $\Omega$. Moreover, Table~\ref{tab:effectivity_index_of_Lshape} presents the values of the effectivity index which show that $I_{\mathrm{eff}}$ for the adaptive method remains bounded and stable---demonstrating that the a posteriori error estimator remains asymptotically proportional to the solution error in energy norm. Figure~\ref{ex2:meshL comparision}(a,b,c) compares the adaptive meshes obtained after 9,8,7 times refinement iterations with $\text{DoF}_h \approx 9000$ for $\mu = 0.1,1,10$, respectively. The meshes illustrate the influence of the surface-diffusion coefficient $\mu$. A larger value of $\mu$ leads to a stronger contribution of the Laplace–Beltrami term, and the corresponding adaptive meshes exhibit a stronger concentration of elements near the boundary. The quasi-optimal convergences of a posteriori error estimators for $\mu = 0.1,1,10$ are displayed in Figure~\ref{ex2:error comparision}.

\begin{table}[htbp]
	\centering
	\caption{Example 2. Effectivity index of the uniform and adaptive mesh refinements.}
	\label{tab:effectivity_index_of_Lshape}
	\begin{tabular}{rcccrccc}
		\toprule
		\multicolumn{4}{c}{Uniform} & \multicolumn{4}{c}{Adaptive} \\
		\cmidrule(lr){1-4} \cmidrule(lr){5-8}
		$\mathrm{DOF}_h$ 
		& $\|\bm{U}-\bm{U}^h\|_{\mathcal{H}^1(\Omega;\Gamma)}$ 
		& $\epsilon_h$ 
		& $I_{\mathrm{eff}}$
		& $\mathrm{DOF}_h$ 
		& $\|\bm{U}-\bm{U}^h\|_{\mathcal{H}^1(\Omega;\Gamma)}$ 
		& $\epsilon_h$ 
		& $I_{\mathrm{eff}}$ \\
		\midrule
		300     & 0.3114 & 4.9662 & 15.95
		& 300    & 0.3114 & 4.9662 & 15.95 \\
		1114    & 0.1700 & 2.5267 & 14.86
		& 1318   & 0.1727 & 1.2393 & 7.18 \\
		4290    & 0.0954 & 1.2806 & 13.42
		& 4538   & 0.0963 & 0.5201 & 5.40 \\
		16834   & 0.0550 & 0.6506 & 11.83
		& 14196  & 0.0480 & 0.2491 & 5.19 \\
		66690   & 0.0325 & 0.3326 & 10.25
		& 21564  & 0.0386 & 0.1906 & 4.94 \\
		265474  & 0.0195 & 0.1718 & 8.79
		& 30526  & 0.0307 & 0.1534 & 5.00 \\
		1059330 & 0.0119 & 0.0901 & 7.55
		& 47502 & 0.0245 & 0.1181 & 4.81 \\
		4232194 & 0.0074 & 0.0481 & 6.54
		& 68206 & 0.0193 & 0.0935 & 4.84 \\
		\bottomrule
	\end{tabular}
\end{table}


\begin{figure}
	\centering
	\includegraphics[scale=0.42]{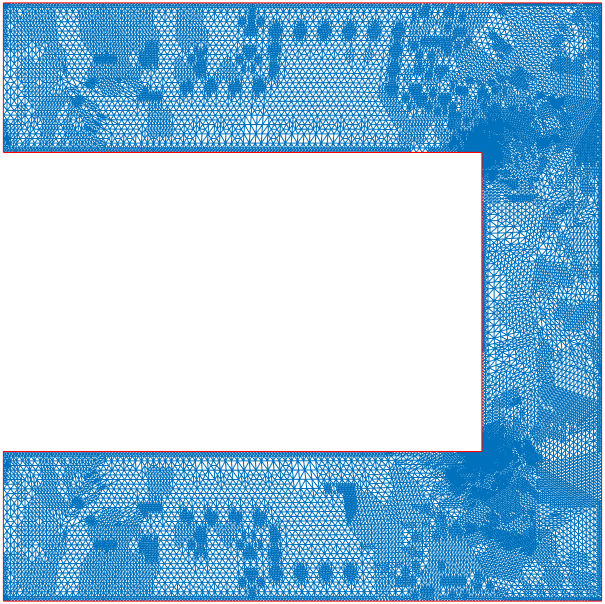} \hspace{0.5cm}
	\includegraphics[scale=0.27]{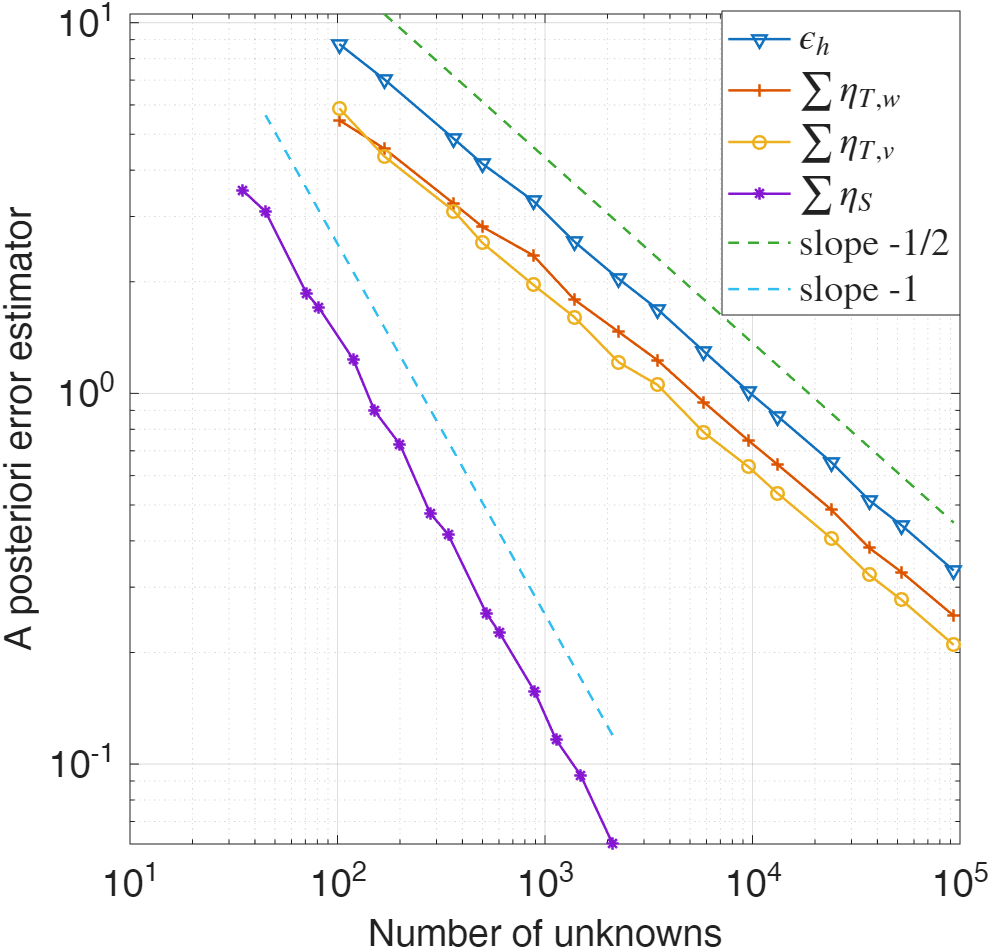} 
	\caption{Example 3. Adaptive mesh (left) and convergence of the a posteriori error estimators (right) for the problem in a U-shaped domain.}
	\label{ex3:result}
\end{figure}

\textbf{Example 3} ({\it U-shaped domain}). Next, we consider the interior transmission problem in a U-shaped domain and choose $k=1, n=8, \mu=2.5-1.5i, \lambda=1.5-2i, \rho_1=1, \rho_2=2$. The boundary data are given by $f=\cos{\frac{\pi x_2}{2}}$ and $g=1+x_1^2+x_2^2$ on $\Gamma$. This case does not provide an analytical solution and due to the re-entrant corners of the U-shaped domain, reduced solution regularity is expected near the corners. Figure~\ref{ex3:result} shows both the mesh after the 12th adaptive refinement and the quasi-optimal convergence of the a posteriori error estimators for the proposed adaptive FEM. 


\textbf{Example 4} ({\it interior transmission eigenvalues}).  In the final example, we study the extension of the proposed adaptive FEM for the numerical calculation of interior transmission eigenvalues for the problem \eqref{eq:Strong Model} with $f=g=0$. A value $k\in \mathbb C$ is called an interior transmission eigenvalue if there exists a nonzero $\bm{U}=(w,v)\in\mathcal H_\Gamma^1(\Omega;\Gamma)$  
such that
\begin{equation*}
	A(\bm{U},\Psi)=0, \qquad \forall\,\Psi=(\psi,\phi) \in\mathcal H_\Gamma^1(\Omega;\Gamma).
\end{equation*}
By separating the $k$ independent and $k$ dependent terms, the sesquilinear form can be re-written as
\begin{equation*}
	A(\bm U,\Psi) = a_e(\bm U,\Psi)-k^2m_e(\bm U,\Psi), 
\end{equation*}
where
\begin{align*}
	a_e(\bm U,\Psi) ={}& \rho_1\int_\Omega \nabla w\cdot\nabla\overline{\psi}\,dx -\rho_2\int_\Omega \nabla v\cdot\nabla\overline{\phi}\,dx + \int_\Gamma \left( \mu\nabla_\Gamma w\cdot \nabla_\Gamma\overline{\psi} + \lambda w\overline{\psi} \right)\,ds,
\end{align*}
and
\begin{equation*}
	m_e(\bm U,\Psi) = \rho_1\int_\Omega n(x)w\overline{\psi}\,dx - \rho_2\int_\Omega v\overline{\phi}\,dx.
\end{equation*}
Then the discrete transmission eigenvalue problem is defined by finding
\begin{equation*}
	(\tau_h,\bm U^h) \in \mathbb  C\times\mathcal H_{\Gamma,h}, \qquad \bm U^h\neq0, 
\end{equation*}
such that
\begin{equation}\label{eq:eigenvalue problem}
	a_e(\bm U^h,\Psi^h) = \tau_h m_e(\bm U^h,\Psi^h), \qquad \forall\,\Psi^h\in\mathcal H_{\Gamma,h},
\end{equation}
where $\tau_h=k_h^2$. Since an eigenfunction is determined only up to a nonzero multiplicative constant, the discrete eigenfunction is normalized as $\|\bm U^h\|_{L^2(\Omega;\Gamma)}=1$.

\begin{figure}
	\centering
	{\includegraphics[width=0.65\textwidth]{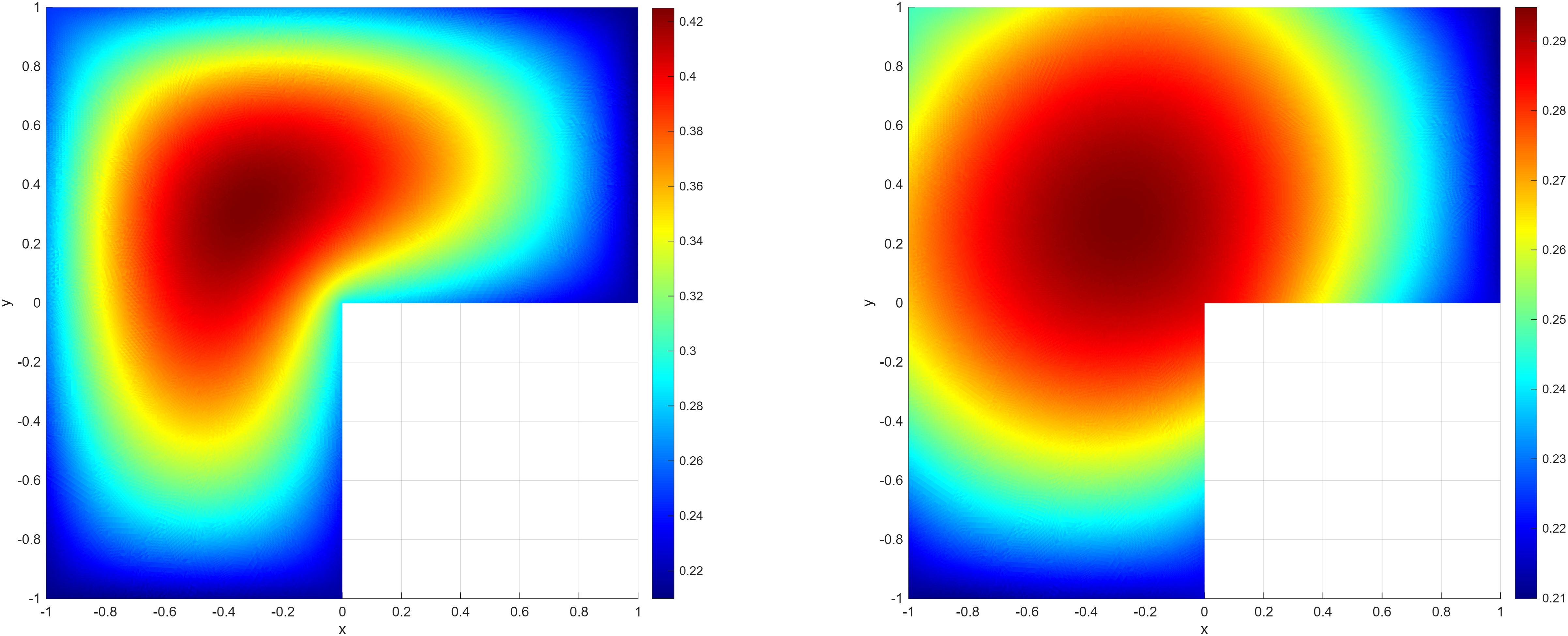}}
	\hspace{0.01 \textwidth}
	{\includegraphics[width=0.272\textwidth]{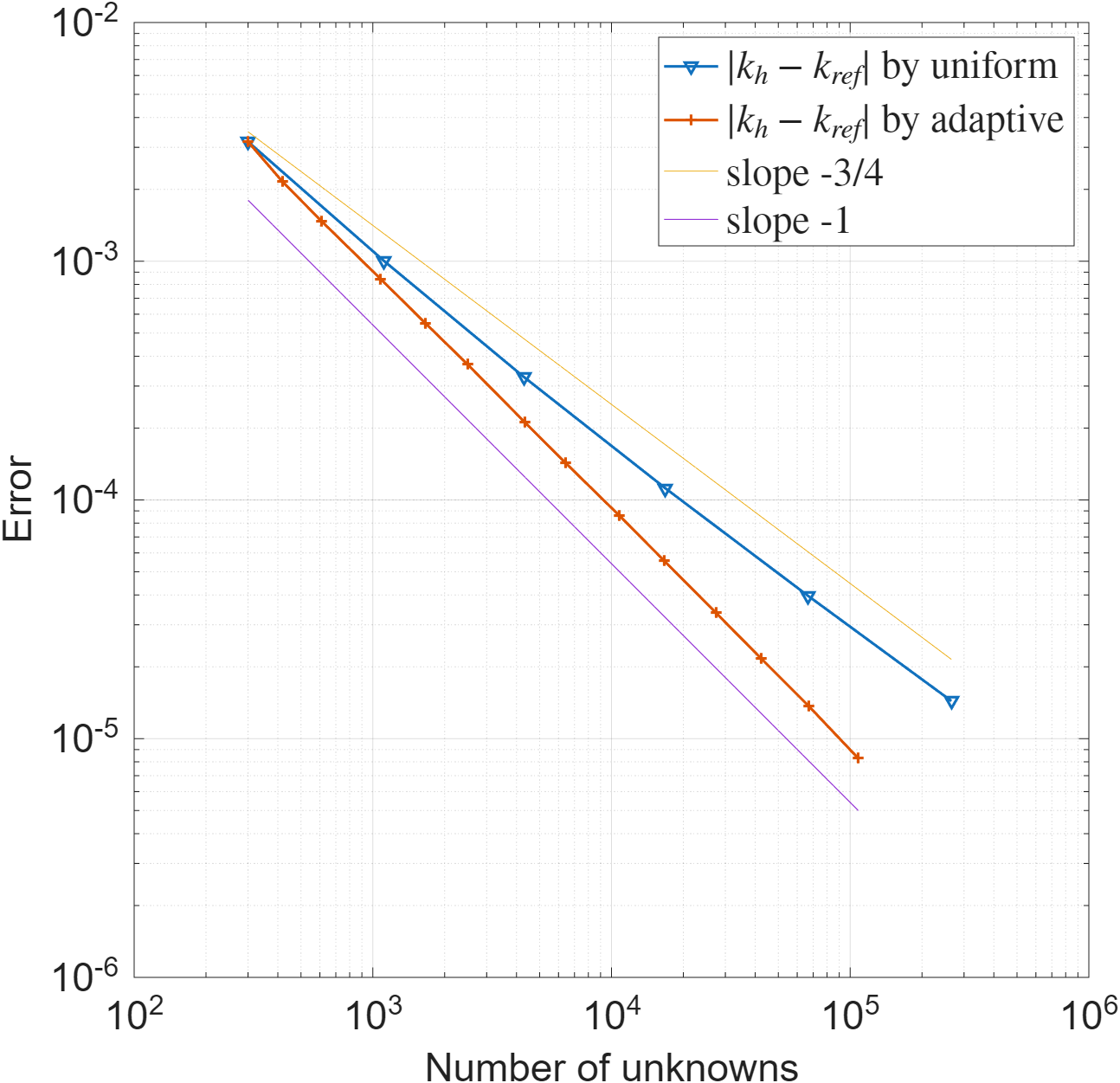}}
	\vspace{-0.5em}
	\caption{Eigenfunctions $w^h$ (left) and $v^h$ (middle) corresponding to eigenvalue $k_{\mathrm{ref}}=0.837617$, and convergence of the errors $E_h$ with respect to $\mathrm {DoF}_h$ for both uniform and adaptive mesh refinements.}
	\label{fig:ex4_error}
\end{figure}

Let $\Omega$ be an L-shaped domain and we set $\rho_1=2,\rho_2=1,n=4,\mu=1.5,\lambda=2.5$. 
The same residual quantities defined in \eqref{eq:Rw}-\eqref{eq:local estimator e} are utilized here for the adaptive computation. The first positive and isolated transmission eigenvalue is considered and tracked under successive mesh refinements. Since no analytical transmission eigenvalue is available for the L-shaped domain, a reference value $k_{\mathrm{ref}}$ is computed on a substantially finer mesh. The eigenvalue error is measured by $E_h = |k_h-k_{\mathrm{ref}}|$ and the convergence rate is calculated through
\begin{equation*}
	\mathrm{rate}= \frac{\log(E_{h-1}/E_{h})}{\log(\mathrm {DoF}_h/\mathrm {DoF}_{h-1})}.
\end{equation*}
The reference eigenvalue $k_{\mathrm{ref}}=0.837617$ is obtained from the 20th adaptive refinement with $\mathrm {DoF}_h=2646 052$. Figure~\ref{fig:ex4_error} shows the eigenfunction corresponding to $k_{\mathrm{ref}}$ and the convergence of the errors $E_h$ with respect to $\mathrm {DoF}_h$ for both uniform and adaptive mesh refinements. Table~\ref{tab:Lshape_TEV} lists the numerical value of the first eigenvalue obtained by the uniform and adaptive finite element approximations. Both results indicate the quasi-optimal convergence $\mathcal{O}\left(\mathrm{DoF}_h^{-1}\right)$ of the adaptive method.

\begin{table}[htbp]
	\centering
	\caption{Example 4. Numerical values of the interior transmission eigenvalue.}
	\begin{tabular}{rcccrccc}
		\toprule
		\multicolumn{4}{c}{Uniform} & \multicolumn{4}{c}{Adaptive} \\
		\cmidrule(lr){1-4} \cmidrule(lr){5-8}
		$\mathrm{DOF}_h$ & $k_h$ & $\left|k_h-k_{\mathrm{ref}}\right|$ & rate 
		& $\mathrm{DOF}_h$ & $k_h$ & $\left|k_h-k_{\mathrm{ref}}\right|$ & rate \\
		\midrule
		300    & 0.840785 & 3.17E-03 & -- 
		& 300    & 0.840785& 3.17E-03 & -- \\
		1114   & 0.838620 & 1.00E-03 & 0.8678 
		& 1076    & 0.838462 & 8.44E-04 & 1.0352 \\
		4290   & 0.837945 & 3.28E-04 & 0.8293 
		& 4332   & 0.837829 & 2.12E-04 & 0.9935 \\
		16834  & 0.837729 & 1.12E-04 & 0.7865 
		& 16578  & 0.837673 & 5.58E-05 & 0.9934 \\
		66690  & 0.837657 & 3.96E-05 & 0.7535 
		& 42388  & 0.837639 & 2.17E-05 & 1.0082 \\
		265474 & 0.837632 & 1.44E-05 & 0.7335 
		& 108102 & 0.837626 & 8.27E-06 & 1.0275 \\
		\bottomrule
	\end{tabular}
	\label{tab:Lshape_TEV}
\end{table}

\section{Conclusion}
\label{sec:5}
This paper discussed the interior transmission problem with GTBCs involving a Laplace–Beltrami term and proposed a new adaptive finite element method with rigorous numerical analysis. Both the a priori error estimates and a posteriori error estimates of the finite element discretization have been derived by combining the $\mathbb{T}$-coercivity technique and finite element analysis. Numerical results demonstrate the efficiency and quasi-optimal convergence of the proposed method as well as its extension to the problem of interior transmission eigenvalues. Developing high-order methods or curved element methods, and the related convergence analysis to treat the high-order derivatives on the boundary more efficiently will be left for future work.


\section*{Funding}
ZZ gratefully acknowledges support from the Postgraduate Scientific Research Innovation Project of Xiangtan University (No. XDCX2026Y272).
LL gratefully acknowledges support from 
the National Natural Science Foundation of China grant (No. 12601798).
WY gratefully acknowledges support from 
the National Natural Science Foundation of China grant (No.12571441) and the 111 Project grant (No. D23017).
TY gratefully acknowledges support from the Beijing Natural Science Foundation (No. JR26005), the National Key R\&D Program of China (No. 2024YFA1016000) and the Strategic Priority Research Program of the Chinese Academy of Sciences (No. XDB0640000).


\end{document}